\documentclass[12pt]{article}
\usepackage{ct}

\newtheorem{res}[theorem]{Result}

\specs{n}{vol (issue)}{year}{}

\dateline{Mar 21, 2025}{Feb 9, 2026}{TBD}



\title{Strongly Regular Graphs with Generalized Denniston and Dual Generalized Denniston Parameters}

\author[1]{Shuxing Li\thanks{Supported by the U.S. National Science Foundation Grant DMS-2452236 and the University of Delaware Research Foundation Strategic Initiative (UDRF-SI) program.}}
\author[2]{James A. Davis}
\author[3]{Sophie Huczynska}
\author[4]{Laura Johnson\thanks{Supported by the Additional Funding Programme for Mathematical Sciences, delivered by EPSRC (EP/V521917/1) and the Heilbronn Institute for Mathematical Research.}}
\author[5]{John Polhill}

\affil[1]{%
Department of Mathematical Sciences, University of Delaware, Newark, DE 19716, USA

\email{shuxingl@udel.edu} (Corresponding author)%
}

\affil[2]{%
Department of Mathematics and Statistics, University of Richmond, Richmond, VA 23173, USA

\email{jdavis@richmond.edu}%
}

\affil[3]{%
School of Mathematics and Statistics, University of St Andrews, North Haugh, St Andrews, Fife, KY16 9SS, Scotland, UK

\email{sh70@st-andrews.ac.uk}%
}

\affil[4]{%
School of Mathematics, University of Bristol, Woodland Road, Bristol, BS8 1UG, England, UK

\email{cc24246@bristol.ac.uk}%
}

\affil[5]{%
Department of Mathematics, Computer Science, and Digital Forensics, Commonwealth University of Pennsylvania, Bloomsburg, PA 17815, USA

\email{jpolhill@commonwealthu.edu}%
}

\newcommand{\Z}{{\mathbb Z}}

\newcommand{\F}{\mathbb{F}}

\newcommand{\al}{\alpha}
\newcommand{\be}{\beta}
\newcommand{\de}{\delta}
\newcommand{\De}{\Delta}
\newcommand{\ga}{\gamma}
\newcommand{\Ga}{\Gamma}
\newcommand{\la}{\lambda}
\newcommand{\om}{\omega}
\newcommand{\tht}{\theta}
\newcommand{\vp}{\varphi}

\newcommand{\Fp}{\F_p}
\newcommand{\Fq}{\F_q}
\newcommand{\Fqm}{\F_{q^m}}
\newcommand{\Fqml}{\F_{q^{m\ell}}}
\newcommand{\Fqmlo}{\F_{q^{m(\ell+1)}}}

\newcommand{\Tr}{\text{Tr}}
\newcommand{\Trqp}{\Tr_{q/p}}
\newcommand{\Trqmp}{\Tr_{q^m/p}}
\newcommand{\Trqmq}{\Tr_{q^m/q}}
\newcommand{\Trqmlp}{\Tr_{q^{m\ell}/p}}
\newcommand{\Trqmlop}{\Tr_{q^{m(\ell+1)}/p}}

\newcommand{\Nqmq}{N_{q^m/q}}
\newcommand{\Nqmlqm}{N_{q^{m\ell}/q^m}}
\newcommand{\Nqmloqm}{N_{q^{m(\ell+1)}/q^m}}

\newcommand{\Cqm}{C^{(e,q^m)}}
\newcommand{\Cqml}{C^{(e,q^{m\ell})}}
\newcommand{\Cqmlo}{C^{(e,q^{m(\ell+1)})}}
\newcommand{\qml}{q^{m\ell}}
\newcommand{\qmlo}{q^{m(\ell+1)}}
\newcommand{\Cl}{C_{\ell}^{\perp}}
\newcommand{\Clo}{C_{\ell+1}^{\perp}}

\newcommand{\pointy}[1]{{\langle#1\rangle}}

\newcommand{\Cay}{\text{Cay}}
\newcommand{\PG}{\text{PG}}
\newcommand{\wt}{\text{wt}}
\newcommand{\es}{\emptyset}
\newcommand{\lan}{\langle}
\newcommand{\ran}{\rangle}
\newcommand{\ol}{\overline}
\newcommand{\pr}{\prime}
\newcommand{\sm}{\setminus}
\newcommand{\wh}{\widehat}
\newcommand{\wti}{\widetilde}
\newcommand{\zn}{\zeta_n}
\newcommand{\zp}{\zeta_p}
\newcommand{\tl}{\tau_{\ell,a,b}}

\begin{document}

\maketitle


\begin{abstract}
We construct two families of strongly regular Cayley graphs, or equivalently, partial difference sets, based on elementary abelian groups. The parameters of these two families are generalizations of the Denniston and the dual Denniston parameters, in contrast to the well known Latin square type and negative Latin square type parameters. The two families unify and subsume a number of existing constructions which have been presented in various contexts such as strongly regular graphs, partial difference sets, projective sets, and projective two-weight codes, notably including Denniston's seminal construction concerning maximal arcs in classical projective planes with even order. Our construction generates further momentum in this area, which recently saw exciting progress on the construction of the analogue of the famous Denniston partial difference sets in odd characteristic.  
\end{abstract}



\section{Introduction}
\label{sec-introduction}

\subsection{Strongly Regular Graphs and Their Relatives}

Strongly regular graphs, according to Peter Cameron, ``lie on the cusp between highly structured and unstructured'' \cite{Cameron}. A $(v,k,\la,\mu)$ {\em strongly regular graph} (SRG) is a simple and undirected $k$-regular graph with $v$ vertices satisfying the following properties
\begin{itemize}
\item[$\bullet$] (edge regular) every two adjacent vertices have $\la$ common neighbours,
\item[$\bullet$] (non-edge regular) every two non-adjacent vertices have $\mu$ common neighbours. 
\end{itemize}
Note that complete or edgeless graphs are trivially strongly regular graphs. Below, we only consider SRGs which are not complete or edgeless. 
 
Let $\Ga$ be a simple and undirected graph with adjacency matrix $A$. The eigenvalues of $\Ga$ are the eigenvalues of the matrix $A$. An eigenvalue of $\Ga$ is \emph{restricted} if it has an eigenvector which is not a multiple of the all-one vector $\bf{1}$. An elegant perspective to analyze SRGs is via their eigenvalues.

\begin{fact}[{\cite[Section 1.1]{BM}}]
\label{fact-SRG}
Let $\Ga$ be a $k$-regular graph with $v$ vertices and adjacency matrix $A$. The following are equivalent:
\begin{itemize}
\item[(1)] $\Ga$ is a $(v,k,\la,\mu)$-SRG.
\item[(2)] $A$ has precisely two restricted eigenvalues, which are the two solutions to the quadratic equation $x^2+(\mu-\la)x+(\mu-k)=0$.
\end{itemize}
\end{fact}

The construction, classification, and nonexistence of SRGs, as well as their connections to other mathematical structures, has been an intensively researched area (see for example the monograph \cite{BM}). A powerful approach to constructing SRGs is to employ Cayley graphs.

\begin{definition}[Cayley graph]
\label{def-Cayley}
Let $G$ be an additively written finite group of order $v$ and $D$ be a subset of $G$ such that the identity $0_G \notin D$ and $D=-D=\{-d \mid d \in D\}$. The Cayley graph on $G$ with connection set $D$, denoted by $\Cay(G,D)$, is the graph with vertices as elements of $G$ and two vertices $g, h \in G$ are adjacent if and only if $g-h \in D$.   
\end{definition}    

\begin{remark}
\label{remark-Cayley}
By Definition~\ref{def-Cayley}, a Cayley graph $\Cay(G,D)$ is simple and undirected. $\Cay(G,D)$ is a $|D|$-regular graph with $G$ being a regular automorphism group of $\Cay(G,D)$. For a subset $S \subset G$, the set of vertices corresponding to $S$ forms a clique in $\Cay(G,D)$ if and only if the set $\{ g-h \mid g, h \in S,  g \ne h\}$ is a subset of $D$. 
\end{remark}

In this paper, we aim to construct two infinite families of strongly regular graphs  whose parameters generalize the Denniston and dual Denniston parameters. For this purpose, we introduce the following concept. 

\begin{definition}[Partial difference set]
\label{def-PDS}
Let $G$ be an additively written group of order $v$ and $D$ be a subset of G with $k$ elements such that the identity $0_G \notin D$. Then $D$ is called a $(v, k,\la,\mu)$ partial difference set (PDS) in $G$ if the expressions $g-h$, for $g$ and $h$ in $D$ with $g \ne h$, represent each element in $D$ exactly $\la$ times and represent each nonidentity element not in $D$ exactly $\mu$ times.
\end{definition}

\begin{remark}
\label{remark-PDS}
\begin{itemize}
\item[(1)] The assumption $0_G \notin D$ in Definition~\ref{def-PDS} is purely technical when considering PDSs: if $D$ is a PDS in $G$, so is $D \cup \{0_G\}$. We insist on $0_G \notin D$, which guarantees the associated Cayley graph $\Cay(G,D)$ does not contain loops.
\item[(2)] In this paper, we only consider $(v,k,\la,\mu)$-PDS $D$ with $\la \ne \mu$, in which $D$ is necessarily fixed by inversion, namely, $D=-D=\{-d \mid d \in D\}$ \cite[Proposition 1.2]{Ma}.  
\item[(3)] Given a $(v,k,\la,\mu)$-PDS $D$ in group $G$, its \emph{complement} $G \sm (D \cup \{ 0_G \})$ is a $(v,v-k-1,v-2-2k+\mu,v-2k+\la)$-PDS in $G$. 
\end{itemize}
\end{remark}

Combining Definition~\ref{def-PDS} and Remark~\ref{remark-PDS}, we have the following equivalence between strongly regular Cayley graphs and PDSs.

\begin{lemma}[{\cite[Proposition 1.1]{Ma}}]
\label{lem-PDS}
A Cayley graph $\Cay(G,D)$ is a $(v,k,\la,\mu)$-SRG if and only if $D$ is a $(v,k,\la,\mu)$-PDS in $G$.
\end{lemma}

Moreover, strongly regular Cayley graphs based on elementary abelian groups are well connected to many other objects, including projective sets and projective two-weight codes. There have been many excellent surveys addressing these connections, including \cite{BM,CK,Ma,MWX}. We also refer to \cite[Section 2]{JL21} for a summary of known constructions of strongly regular Cayley graphs based on finite abelian groups. Next, we summarize the relation between strongly regular Cayley graphs based on elementary abelian groups and their relatives. For standard terminologies concerning projective geometry and coding theory, please refer to \cite{CK,Ma}. We first introduce the definitions of projective sets and projective two-weight codes. From now on, we always use $q$ to represent a prime power.

Let $\PG(m-1,q)$ be the $(m-1)$-dimensional projective space over $\Fq$. The point set of $\PG(m-1,q)$ is $\{ \lan y_i \ran \mid 1 \le i \le \frac{q^m-1}{q-1} \}$, where each $\lan y_i \ran$ is an one-dimensional vector subspace of $\Fq^m$ over $\Fq$ with $\{ y_i \}$ as a basis, and $\{ \lan y_i \ran \mid 1 \le i \le \frac{q^m-1}{q-1} \}$ consists of all one-dimensional vector subspaces of $\Fq^m$ over $\Fq$. A hyperplane in $\PG(m-1,q)$ corresponds to a $(m-1)$-dimensional vector subspace of $\Fq^m$ over $\Fq$. The set of all hyperplanes in $\PG(m-1,q)$ corresponds to the set of all $(m-1)$-dimensional vector subspaces of $\Fq^m$ over $\Fq$.

\begin{definition}[Projective set]
A projective $(n, m, h_1,h_2)$ set $O$ is a proper, non-empty set of $n$ points of the projective space $\PG(m-1,q)$ with the property that each hyperplane of $\PG(m-1,q)$ meets $O$ in either $h_1$ or $h_2$ points.
\end{definition}

For $0 \le m \le n$, let $C$ be an $m$-dimensional vector subspace of $\Fq^n$ over $\Fq$. Then $C$ is an $[n,m]_q$ code. For $c=(c_1,c_2,\ldots,c_n) \in C$, define the \emph{weight} of $c$ to be 
$$
\wt(c)=|\{ 1 \le i \le n \mid c_i \ne 0 \}|.
$$
$C$ is called a \emph{two-weight code} if all the vectors in $C$ have precisely two distinct nonzero weights.
\begin{definition}[Projective two-weight code]
For $1 \le i \le n$, let $y_i \in \Fq^m$. The $[n,m]_q$ code $C$ defined by
$$
C=\{ (x \cdot y_1, x \cdot y_2, \ldots, x \cdot y_n) \mid x \in \Fq^m\},
$$
where $\cdot$ is the usual dot product in $\Fq^m$, is a projective two-weight code if both of the following hold:
\begin{itemize}
\item[(1)] For every $1 \le i, j \le n$ with $i \ne j$, $y_i$ and $y_j$ are linearly independent over $\Fq$, namely, $C$ is a projective code.
\item[(2)] $C$ is a two-weight code. 
\end{itemize}
\end{definition}

Now we are ready to describe the connection between strongly regular Cayley graphs based on elementary abelian groups and their relatives.

\begin{lemma}[{\cite[Section 3]{CK}}]
\label{lem-connection}
Let $q$ be a prime power. Let $G=\Fq^m$. Let $O=\{ \lan y_i \ran \mid 1 \le i \le n\}$ be a set of distinct points in $\PG(m-1,q)$. Let $D=\{ \al y_i \mid 1 \le i \le n,  \al \in \Fq^* \}$. Then the following are equivalent:
\begin{itemize}
\item[(1)] $O$ is a projective $(n,m,n-w_1,n-w_2)$ set in $\PG(m-1,q)$.
\item[(2)] The $[n,m]_q$ code $C=\{ (x \cdot y_1, x \cdot y_2, \ldots, x \cdot y_n) \mid x \in \Fq^m\}$ is a projective two-weight code with nonzero weights $w_1$ and $w_2$. 
\item[(3)] $\Cay(G,D)$ is a $(v,k,\la,\mu)$-SRG with
\begin{equation}
\label{eqn-para}
\begin{split}
v=&q^m, \\
k=&n(q-1), \\
\la=& k^2+3k-q(w_1+w_2)-kq(w_1+w_2)+q^2w_1w_2, \\
\mu=& k^2+k-kq(w_1+w_2)+q^2w_1w_2.
\end{split}
\end{equation}
\item[(4)] $D$ is a $(v,k,\la,\mu)$-PDS in $G$, where $(v,k,\la,\mu)$ satisfies equations in (\ref{eqn-para}).
\end{itemize}
\end{lemma}

Below, we will focus on a particular family of SRGs with so-called Denniston parameters and describe the state-of-the-art. 

\subsection{Strongly Regular Cayley Graphs with Denniston Parameters and Their Extensions}

In view of Lemma~\ref{lem-PDS}, the parameters of strongly regular Cayley graphs exactly mirror those of the corresponding PDS. Therefore, in this subsection, we do not distinguish strongly regular Cayley graphs and PDSs, and omit mentioning them explicitly when describing their parameter families. 

Most known strongly regular Cayley graphs, or equivalently, PDSs, fall into the following two well known families. 
\begin{itemize}
\item[$\bullet$] Latin square type family with parameters $(n^2,r(n-1),n+r^2-3r,r^2-r)$ for $n \ge 1$ and $r \ge 0$.  
\item[$\bullet$] Negative Latin square type family with parameters $(n^2,r(n+1),-n+r^2+3r,r^2+r)$ for $n \ge 1$ and $r \ge 0$.  
\end{itemize}
We note that in both of the above families, the order of the group is necessarily a square. Given that, the Denniston family stands out for its distinct parameters as compared with the Latin square type and the negative Latin square type families, as well as for its intrinsic connection with finite geometry. Actually, Denniston's pioneering construction concerns maximal arcs in classical projective plane of even order. Specifically, let $n \ge 2$ and $d \ge 1$. An \emph{$(n,d)$-arc} in the classical projective plane $\PG(2,q)$ is a set of $n$ points, of which no $d+1$ points are collinear. For an $(n,d)$-arc in $\PG(2,q)$, we have 
$$ 
n \le 1+(q+1)(d-1),
$$
where the $(n,d)$-arc is \emph{maximal} if $n=1+(q+1)(d-1)$. Note that each line in $\PG(2,q)$ intersects a maximal $(n,d)$-arc at either $0$ or $d$ points. Then a maximal $(n,d)$-arc is a projective $(n,3,0,d)$ set in $\PG(2,q)$. Considering \emph{nontrivial} maximal $(n,d)$-arcs where $1 < d < q$, it can be shown that $d \mid q$. In the classical projective plane $\PG(2,q)$ with $q=2^m$, for each $1 \le r \le m-1$, Denniston constructed maximal $(n,2^r)$-arcs \cite{Denniston}, which led to the parameter family $(2^{3m}, (2^{m+r} - 2^m + 2^r)(2^m-1), 2^m-2^r+(2^{m+r}-2^m+2^r)(2^r-2), (2^{m+r}-2^m+2^r)(2^r-1))$  in $\F_2^{3m}$.   

Following Denniston's construction, a natural question is whether strongly regular Cayley graphs and PDSs with Denniston parameters can exist in the additive group of a finite field with odd characteristic. In short, we say the question concerns Denniston parameters in odd characteristic. Notably, a famous result due to Ball, Blokhuis, and Mazzocca establishes the nonexistence of nontrivial maximal arcs in classical projective planes of odd order \cite{BBM}. For a long time, this result seemed to rule-out the possibility of Denniston parameters in odd characteristic. Consequently, the recent discovery of SRGs and PDSs with Denniston parameters in odd characteristic, simultaneously and independently by two groups of researchers, came as a huge surprise. Specifically, Davis, Huczynska, Johnson, and Polhill employed cyclotomic classes and quadratic forms over finite fields to construct the parameter family $(p^{3m}, (p^{m+r} - p^m + p^r)(p^m-1), p^m-p^r+(p^{m+r}-p^m+p^r)(p^r-2), (p^{m+r}-p^m+p^r)(p^r-1))$ in $\F_{p^m} \times \F_{p^{2m}}$, where $p$ is an odd prime, $m \ge 2$, and $r \in \{ 1, m-1 \}$ \cite[Theorems 5.1, 5.2]{DHJP}. De Winter utilized both geometric and algebraic approaches to construct projective sets in $\PG(3m-1,q)$, which amount to parameter family $(q^{3m},(q^{m+1}-q^m+q)(q^m-1),q^m-q+(q^{m+1}-q^m+q)(q-2),(q^{m+1}-q^m+q)(q-1))$ in $\Fqm \times \F_{q^{2m}}$, where $q$ is a prime power and $m \ge 2$ \cite[Theorem 2]{DeWinter}. Recently, Bao, Xiang, and Zhao extended both constructions and presented the parameter family $(q^{3m}, (q^{m+r} - q^m + q^r)(q^m-1), q^m-q^r+(q^{m+r}-q^m+q^r)(q^r-2), (q^{m+r}-q^m+q^r)(q^r-1))$ in $\Fqm \times \F_{q^{2m}}$, where $q$ is a prime power, $m \ge 2$, and $1 \le r \le m-1$ \cite[Theorem 3.2]{BXZ}.  Given that $q$ could be an arbitrary prime power and $r$ could vary over the full possible range $1 \le r \le m-1$, this construction established the existence of all possible Denniston parameter families in elementary abelian groups. Moreover, a neat explanation was presented in \cite[Section 1]{BXZ}, illustrating that maximal arcs in projective planes of odd order are stronger than SRGs and PDSs with Denniston parameters in odd characteristic. Therefore, the Denniston parameter family in odd characteristic does not contradict the nonexistence result in \cite{BBM}.

The multiple perspectives from which SRGs and their relatives may be viewed, showcased in Lemma~\ref{lem-connection}, speaks for the charm of the topic. On the other hand, it is not uncommon that some results described via one perspective have not been recognized by others. For instance, as De Winter observed \cite[Section 1]{DeWinter}, in a less well known paper \cite{BE}, Bierbrauer and Edel effectively obtained the same parameter family as that of \cite{BXZ}, via an intricate coding-theoretic approach. During the preparation of this paper, the authors came across another less noticed paper by Ott \cite{OttDM}, which also extended Denniston's construction in a different direction. While Denniston's construction was based on group $\F_{2^{3m}}=\F_{2^m} \times \F_{2^{2m}}$, Ott considered a much more general construction based on group $\Fqml \times \Fqmlo$, where $q$ is a prime power and $m, \ell \ge 1$. As we shall see in Section~\ref{sec-statement}, our construction of generalized Denniston and dual generalized Denniston parameter families further extends Ott's construction.

\subsection{SRGs and PDSs with Generalized Denniston and Dual Generalized Denniston Parameters}

We first describe the family of SRGs and PDSs with generalized Denniston parameters.

\begin{theorem}
\label{theorem-main}
Let $s \ge 1$ and $q=p^s$ be a power of prime $p$. Let $m$ and $\ell$ be positive integers. For each $0 \leq r \leq m$, there exists a $(v,k_r,\la_r,\mu_r)$-PDS $D_r$ in the elementary abelian $p$-group $G=\Z_p^{sm(2\ell+1)}$, where
\begin{equation*}
\begin{split}
v = & q^{m(2 \ell+1)}, \\
k_r  = & \frac{(q^r-1)(q^{m\ell}-1)(q^{m(\ell+1)}-1)}{q^m-1} + q^{m\ell}-1, \\
\la_r = & \frac{(q^r-1)(q^{m(\ell+1)}-1)}{q^m-1}\Big(\frac{(q^r-1)(q^{m\ell}-1)}{q^m-1}-1\Big) + q^{m\ell}-2, \\
\mu_r = & \frac{(q^r-1)(q^{m\ell}-1)}{q^m-1}\Big(\frac{(q^r-1)(q^{m(\ell+1)}-1)}{q^m-1}+1\Big). 
\end{split}
\end{equation*}
Consequently, there exists a $(v,k_r,\la_r,\mu_r)$-SRG $\Cay(G,D_r)$.
\end{theorem}

\begin{remark}
\label{remark-main}
The infinite family of SRGs and PDSs in Theorem~\ref{theorem-main} have the generalized Denniston parameters, as they extend a series of previous constructions initiated by Denniston:
\begin{itemize}
\item[(1)] Applying $p=2$, $s=1$, $\ell=1$, and $1 \le r \le m-1$ to Theorem~\ref{theorem-main}, we obtain SRGs and PDSs with the same parameters as Denniston's initial work \cite{Denniston}.
\item[(2)] Applying $m=2$ and $r=1$ to Theorem~\ref{theorem-main}, we recover the construction by Momihara \cite[Equation (1.1)]{M}, which is of the negative Latin square type. We note that Momihara's construction \cite[Proposition 2]{M} extends that of Polhill \cite[Corollaries 5.1, 5.2]{P}. 
\item[(3)] Applying $r=1$ to Theorem~\ref{theorem-main}, we recover the construction by Ott \cite[Theorem 3]{OttDM}. We note that Ott's construction extends that of Fern\'{a}ndez-Alcober, Kwahira, and  Mart\'{\i}nez \cite[Proposition 4.3]{FKM}.
\item[(4)] Applying $p$ being an odd prime, $s=1$, $\ell=1$, and $r \in \{1,m-1\}$ to Theorem~\ref{theorem-main}, we recover the construction by Davis, Huczynska, Johnson, and Polhill \cite[Theorem 5.1]{DHJP}.
\item[(5)] Applying $\ell=1$ and $1 \le r \le m-1$ to Theorem~\ref{theorem-main}, we recover the construction by Bao, Xiang, and Zhao \cite[Theorem 3.2]{BXZ}.
\end{itemize}
\end{remark}

Given a strongly regular Cayley graph $\Cay(G,D)$ based on a finite abelian group $G$, or equivalently, a PDS $D$ in a finite abelian group $G$, there exists a dual strongly regular Cayley graph $\Cay(\wh{G},D^+)$, or equivalently, a dual PDS $D^+$ in the character group $\wh{G}$ of $G$. Below, we describe the family of SRGs and PDSs with dual generalized Denniston parameters.

\begin{theorem}
\label{theorem-dualmain}
Let $s \ge 1$ and $q=p^s$ be a power of prime $p$. Let $m$ and $\ell$ be positive integers. For each $0 \leq r \leq m$, there exists a $(v^+,k_r^+,\la_r^+,\mu_r^+)$-PDS $D_r^+$ in the elementary abelian $p$-group $G=\Z_p^{sm(2\ell+1)}$, where
\begin{equation*}
\begin{split}
v^+ = & q^{m(2 \ell+1)}, \\
k_r^+  = & \frac{(q^m-q^{m-r})(q^{m\ell}-1)(q^{m(\ell+1)}-1)}{q^m-1} + q^{m(\ell+1)}-1, \\
\lambda_r^+ = & \frac{(q^m-q^{m-r})(q^{m(\ell+1)}-1)(q^{m(\ell+1)}-q^{m(\ell+1)-r}+q^{m-r}-1)}{(q^m-1)^2} \\
                        & + \frac{q^{m(\ell+2)-r}+q^{m(\ell+1)-r}-2q^{m(\ell+1)}-2q^{m-r}+2}{q^m-1}, \\
\mu_r^+ = & \frac{(q^m-q^{m-r})(q^{m(\ell+1)}-1)(q^{m(\ell+1)}-q^{m(\ell+1)-r}+q^{m-r}-1)}{(q^m-1)^2}.
\end{split}
\end{equation*}
Consequently, there exists a $(v^+,k_r^+,\la_r^+,\mu_r^+)$-SRG $\Cay(G,D_r^+)$.
\end{theorem}
Note that the group $G=\Z_p^{sm(2\ell+1)}$ in Theorem~\ref{theorem-main} has character group $\wh{G}$, which is again isomorphic to $\Z_p^{sm(2\ell+1)}$. That is why we phrase Theorem~\ref{theorem-dualmain} with respect to the group $\Z_p^{sm(2\ell+1)}$. The SRGs and PDSs in Theorem~\ref{theorem-dualmain} have the dual generalized Denniston parameters, which are in general distinct from the generalized Denniston parameters. The standard procedure of deriving the dual PDS $D^+$ from the original PDS $D$ will be detailed in Lemma~\ref{lem-dual}.

\begin{remark}
\label{remark-dualmain}
As we shall see in Theorem~\ref{theorem-dualmain2}, $D_r^+$ can be alternatively constructed as the complement of some (properly constructed) $D_{m-r}$. In this sense, the dual generalized Denniston family in Theorem~\ref{theorem-dualmain} is not that new. However, for completeness, we include this family. We also note the corresponding SRG $\Cay(G,D_r^+)$ has an interesting property concerning its maximum clique, detailed in Theorem~\ref{theorem-dualmain2}. 
\end{remark}

\begin{remark}
\label{remark-relation}
In general, SRGs and PDSs with generalized Denniston and dual generalized Denniston parameters are not of the Latin square or negative Latin square type, say, for instance, when $m$ is odd. On the other hand, when $m$ is even and $r=\frac{m}{2}$, the generalized Denniston family satisfies 
$$
k_{\frac{m}{2}}=\frac{(q^{m\ell}-1)(q^m-q^{\frac{m}{2}})(q^{\frac{m}{2}(2\ell+1)}+1)}{q^m-1}
$$ 
and the dual generalized Denniston family satisfies
$$
k_{\frac{m}{2}}^+=\frac{(q^{m(\ell+1)}-1)(q^{\frac{m}{2}}-1)(q^{\frac{m}{2}(2\ell+1)}+1)}{q^m-1}.
$$
Therefore, the corresponding SRGs and PDSs are of the negative Latin square type.
\end{remark}

In view of Lemma~\ref{lem-connection}, we have the following corollary regarding projective sets derived from Theorems~\ref{theorem-main} and~\ref{theorem-dualmain}.

\begin{corollary}
\label{corollary-projset}
Let $q$ be a prime power. Let $m$ and $\ell$ be positive integers. For each $0 \leq r \leq m$, the following exist.
\begin{itemize}
\item[(1)] A projective $(n_r,m(2\ell+1),h_{r,1},h_{r,2})$ set in $\PG(m(2\ell+1)-1,q)$, where
\begin{equation*}
\begin{split}
n_r  = & \frac{(q^r-1)(q^{m\ell}-1)(q^{m(\ell+1)}-1)}{(q-1)(q^m-1)} + \frac{q^{m\ell}-1}{q-1}, \\
h_{r,1} = & \frac{(q^{m\ell-1}-1)(q^{m(\ell+1)+r}-q^{m(\ell+1)}+q^m-q^r)}{(q-1)(q^m-1)}, \\
h_{r,2} = & \frac{(q^{m\ell}-1)(q^{m(\ell+1)+r-1}-q^{m(\ell+1)-1}+q^m-q^r)}{(q-1)(q^m-1)}. 
\end{split}
\end{equation*}
\item[(2)] A projective $(n_r^+,m(2\ell+1),h_{r,1}^+,h_{r,2}^+)$ set in $\PG(m(2\ell+1)-1,q)$, where
\begin{equation*}
\begin{split}
n_r^+  = & \frac{(q^m-q^{m-r})(q^{m\ell}-1)(q^{m(\ell+1)}-1)}{(q-1)(q^m-1)} + \frac{q^{m(\ell+1)}-1}{q-1}, \\
h_{r,1}^+ = & \frac{(q^{m(\ell+1)-1}-1)(q^{m(\ell+1)}-q^{m(\ell+1)-r}+q^{m-r}-1)}{(q-1)(q^m-1)}, \\
h_{r,2}^+ = & \frac{(q^{m(\ell+1)}-1)(q^{m(\ell+1)-1}-q^{m(\ell+1)-r-1}+q^{m-r}-1)}{(q-1)(q^m-1)}.
\end{split}
\end{equation*}
\end{itemize}
\end{corollary}

\begin{remark}
\label{remark-projset}
Applying $\ell=1$ and $r=1$ to Corollary~\ref{corollary-projset}(1), we recover the projective $\big(\frac{(q^m-1)(q^{m+1}-q^m+q)}{q-1},3m,\frac{q^{2m}-q^{2m-1}-q^{m+1}+2q^m-q}{q-1},\frac{q^{2m}-q^{2m-1}+q^m-q}{q-1}\big)$ set constructed by De Winter \cite[Theorem 4.1]{DeWinter}.
\end{remark}

An intriguing feature of De Winter's construction is that the obtained projective set has parameters corresponding to the so-called blow up of a hypothetical maximal $(n,q)$-arc in $\PG(2,q^m)$. Note that a point in $\PG(2,q^m)$, being a one-dimensional vector space over $\Fqm$, corresponds to an $m$-dimensional vector space over $\Fq$, which again corresponds to an $(m-1)$-space in $\PG(3m-1,q)$. Following this manner, the blow up process transforms objects in $\PG(2,q^m)$ to objects in $\PG(3m-1,q)$, where a point set in $\PG(2,q^m)$ is transformed to a set of distinct $(m-1)$-spaces in $\PG(3m-1,q)$ and a line in $\PG(2,q^m)$ is transformed to a $(2m-1)$-space in $\PG(3m-1,q)$. 

We now illustrate that a projective  
$$
\Big( \frac{(q^m-1)(q^{m+1}-q^m+q)}{q-1}, 3m, \frac{q^{2m}-q^{2m-1}-q^{m+1}+2q^m-q}{q-1}, \frac{q^{2m}-q^{2m-1}+q^m-q}{q-1} \Big)
$$
set can be obtained by blowing up a hypothetical maximal $(n,q)$-arc in $\PG(2,q^m)$. A hypothetical maximal $(n,q)$-arc consists of $n=1+(q^m+1)(q-1)=q^{m+1}-q^m+q$ points, so that every line in $\PG(2,q^m)$ intersects the maximal $(n,q)$-arc in either $0$ or $q$ points. Hence, the blow up of the maximal $(n,q)$-arc leads to a point set in $\PG(3m-1,q)$ with size
$$
(q^{m+1}-q^m+q)\frac{q^m-1}{q-1}=\frac{(q^m-1)(q^{m+1}-q^m+q)}{q-1}.
$$
And this point set intersects the hyperplanes in $\PG(3m-1,q)$ with exactly two sizes: 
\begin{align*}
(q^{m+1}-q^m+q)\frac{q^{m-1}-1}{q-1}&=\frac{q^{2m}-q^{2m-1}-q^{m+1}+2q^m-q}{q-1}, \\
(q^{m+1}-q^m)\frac{q^{m-1}-1}{q-1}+q\frac{q^m-1}{q-1}&=\frac{q^{2m}-q^{2m-1}+q^m-q}{q-1}.
\end{align*}

With the above preparation, we claim that Corollary~\ref{corollary-projset}(1) generalizes De Winter's construction in the following sense. Applying $\ell=1$ and $1 \le r \le m-1$ to Corollary~\ref{corollary-projset}(1), we obtain a projective
$$
\Big(\frac{(q^m-1)(q^{m+r}-q^m+q^r)}{q-1},3m, \frac{(q^{m-1}-1)(q^{m+r}-q^m+q^r)}{q-1}, \frac{q^{2m+r-1}-q^{2m-1}+q^m-q^r}{q-1}\Big)
$$
set in $\PG(3m-1,q)$. Following the same argument as above, this projective set has parameters corresponding to the blow up of a hypothetical maximal $(n,q^r)$-arc in $\PG(2,q^m)$. We shall note that, as described later in Remark~\ref{remark-main2}(5), the projective set in Corollary~\ref{corollary-projset}(1) with $1 \le r \le m-1$ cannot be the blow up of a point set in $\PG(2,q^m)$.

Remarkably, De Clerck, Delanote, Hamilton, and Mathon discovered a projective $(84,6,21,30)$ set in $\PG(5,3)$ through their careful study of the \emph{perp-system} \cite{DDHM}. Their example is formed as the union of $21$ joint lines of $\PG(5,3)$ and exhibits many intriguing properties \cite[Example~2]{DDHM}. Consequently, this example is genuinely distinct from the projective $(84,6,21,30)$ set arising from Corollary~\ref{corollary-projset}(1) with $q=3$, $m=2$, $\ell=1$, and $r=1$.

Finally, in view of Lemma~\ref{lem-connection}, we have the following corollary regarding projective two-weight codes derived from Theorems~\ref{theorem-main} and~\ref{theorem-dualmain}.

\begin{corollary}
\label{corollary-code}
Let $q$ be a prime power. Let $m$ and $\ell$ be positive integers. For each $0 \leq r \leq m$, the following exist.
\begin{itemize}
\item[(1)] A projective two-weight $[n_r,m(2\ell+1)]_q$ code with two weights $w_{r,1}$ and $w_{r,2}$, where
\begin{equation*}
\begin{split}
n_r  = & \frac{(q^r-1)(q^{m\ell}-1)(q^{m(\ell+1)}-1)}{(q-1)(q^m-1)} + \frac{q^{m\ell}-1}{q-1}, \\
w_{r,1} = & \frac{q^{m\ell-1}(q^{m(\ell+1)+r}-q^{m(\ell+1)}+q^m-q^r)}{q^m-1}, \\
w_{r,2} = & \frac{q^{m(\ell+1)-1}(q^r-1)(q^{m\ell}-1)}{q^m-1}. 
\end{split}
\end{equation*}
\item[(2)] A projective two-weight $[n_r^+,m(2\ell+1)]_q$ code with two weights $w_{r,1}^+$ and $w_{r,2}^+$, where
\begin{equation*}
\begin{split}
n_r^+  = & \frac{(q^m-q^{m-r})(q^{m\ell}-1)(q^{m(\ell+1)}-1)}{(q-1)(q^m-1)} + \frac{q^{m(\ell+1)}-1}{q-1}, \\
w_{r,1}^+ = & \frac{q^{m(\ell+1)-r-1}(q^{m(\ell+1)+r}-q^{m(\ell+1)}+q^m-q^r)}{q^m-1}, \\
w_{r,2}^+ = & \frac{(q^r-1)(q^{m(\ell+1)}-1)(q^{m(\ell+1)-r}-q^{m(\ell+1)-r-1})}{(q-1)(q^m-1)}. 
\end{split}
\end{equation*}
\end{itemize}
\end{corollary}

\begin{remark}
\label{remark-code}
Applying $\ell=1$ to Corollary~\ref{corollary-code}(1), we obtain a projective two-weight 
$$
[\frac{q^r(q^m-1)(q^m-q^{m-r}+1)}{q-1},3m]_q
$$ 
code with weights $q^{m+r-1}(q^m-q^{m-r}+1)$ and $q^{2m+r-1}-q^{2m-1}$, which has the same parameter as the projective two-weight code constructed by Bierbrauer and Edel \cite{BE}.
\end{remark}

The rest of the paper is organized as follows. Section~\ref{sec-prelim} provides necessary background, including cyclotomic classes, additive and multiplicative characters over finite fields, Gauss sums, and dual PDSs. Section~\ref{sec-statement} presents explicit constructions of SRGs and PDSs with generalized Denniston and dual generalized Denniston parameters. The explicit expressions for the generalized Denniston and the dual generalized Denniston families help in articulating the connection between these two new families and the existing constructions, resulting in an affirmative answer to an open problem proposed by De Winter \cite{DeWinter}. Section~\ref{sec-proof} contains the proof of our main result, establishing Theorems~\ref{theorem-main} and \ref{theorem-dualmain}. The proof relies on a careful character-theoretic analysis, inspired by the enlightening work of Momihara and Xiang \cite{MX}. Section~\ref{sec-conclusion} contains some concluding remarks and several open problems.

\section{Preliminaries}
\label{sec-prelim}

In this section, we first give a brief introduction to characters and character sums, which are instrumental in showing that a subset of a finite abelian group is a  PDS. We then describe the standard process of deriving the dual PDS from the original one. As our main construction relies on a careful combination of cyclotomic classes of finite fields, we subsequently define cyclotomic classes as well as additive and multiplicative characters over finite fields. When applying additive characters to a union of cyclotomic classes, Gauss sums naturally come into play. Some basic facts about Gauss sums are also mentioned. 

\subsection{Characters and Character Sums}

Let $G$ be a finite abelian group with identity $1_G$. The \emph{exponent} of a group $G$ is the smallest positive integer~$n$ such that $g^n = 1_G$ for each $g \in G$.  A \emph{character} $\chi$ of a group $G$ is a group homomorphism from $G$ to the multiplicative group~$\pointy{\zn}$ generated by a primitive $n$-th root of unity $\zn$.  

Denote the set consisting of all characters on $G$ by $\wh{G}$. For $\chi_1, \chi_2 \in \wh{G}$, define an operation $\circ$ such that
$$
(\chi_1 \circ \chi_2)(g)=\chi_1(g)\chi_2(g), \forall g \in G. 
$$
Then $(\wh{G}, \circ)$ forms a group isomorphic to $G$, called the \emph{character group} of $G$. The character group $\wh{G}$ is isomorphic to $G$. Within $\wh{G}$, the \emph{principal character} is the character that maps all elements of $G$ to $1$; all the other characters of $\wh{G}$ are called \emph{nonprincipal characters}. For $\chi \in \wh{G}$ and a positive integer $t$, we use $\chi^t$ to represent the $t$-fold multiplication under the operation $\circ$, namely, 
$$
\chi^{t}=\underbrace{\chi\circ\chi\circ\cdots\circ\chi}_{\text{$t$ occurrences of $\chi$}}.
$$
Moreover, $\chi^0$ is the principal character of $\wh{G}$. We use $\chi^{-1}$ to denote the inverse of $\chi$. Consequently, for each $g \in G$, $\chi^{-1}(g)=\ol{\chi(g)}$. Let $H$ be a subgroup of $G$. A character $\chi$ is \emph{principal on $H$} if $\chi$ maps each element of $H$ to $1$. Define $H^{\perp}=\{ \chi \in \wh{G} \mid \mbox{$\chi$ is principal on $H$}\}$.

Let $\chi \in \wh{G}$ and $D$ a subset of $G$. Define the \emph{character sum} of $D$ with respect to $\chi$ as $\chi(D)=\sum_{d \in D} \chi(d)$. For instance, given a subset $D \subset G$, for the principal character $\chi_0 \in \wh{G}$, the character sum of $D$ with respect to $\chi_0$ is $\chi_0(D)=|D|$. The following lemma provides a neat characterization of PDSs via their character sums.

\begin{lemma}
[{\cite[Corollary.~3.3]{Ma}}]
\label{lem-pdschar}
Let $G$ be an abelian group of order $v$. Let $D$ be a $k$-subset of~$G$ for which $1_G \notin D$.
Let $\la,\mu$ be nonnegative integers satisfying
$k^2=k-\mu+(\la-\mu)k+\mu v$ and $(\la-\mu)^2+4(k-\mu) \ge 0$.
Then $D$ is a $(v,k,\la,\mu)$-PDS in $G$ if and only if
$$
\chi(D) = \frac{\la-\mu \pm \sqrt{(\la-\mu)^2+4(k-\mu)}}{2}  \mbox{for all nonprincipal characters $\chi$ of $G$}.
$$
\end{lemma}
Lemma~\ref{lem-pdschar} indicates that a subset $D$ of a finite abelian group $G$ is a PDS if and only if the character sums $\{\chi(D) \mid \chi \in \wh{G}, \chi \mbox{\;nonprincipal}\}$ have exactly two distinct values. This result is critical for our subsequent analysis in Section~\ref{sec-proof}. For an expository introduction to characters and character sums, please refer to \cite[Section 2]{JL24}. 

The following lemma due to Delsarte indicates that a PDS in finite abelian group $G$ automatically generates a dual PDS in the character group $\wh{G}$.

\begin{lemma}[{\cite[Theorem 3.4]{Ma}}]
\label{lem-dual}
Let $G$ be an abelian group of order $v$ and $D$ be a $(v,k,\la,\mu)$-PDS in $G$ with $D \ne \es$ and $D \ne G$. Set $\be=\la-\mu$, $\ga=k-\mu$, and $\De=\be^2+4\ga$. Then the dual $D^+$ of $D$, defined as
$$
D^+=\{ \chi \in \wh{G} \mid \mbox{$\chi$ is nonprincipal and $\chi(D)=\tfrac{1}{2}(\la-\mu + \sqrt{(\la-\mu)^2+4(k-\mu)})$} \}
$$
is a $(v^+,k^+,\la^+,\mu^+)$-PDS in $\wh{G}$ with
\begin{align*}
v^+=&v, \\
k^+=& \frac{(\sqrt{\De}-\be)(v-1)-2k}{2\sqrt{\De}}, \\
\la^+=&k^{+}+\frac{(v-2k+\be-\sqrt{\De})^2-v^2}{4\De}+\frac{v-2k+\be-\sqrt{\De}}{\sqrt{\De}}, \\
\mu^+=&k^{+}+\frac{(v-2k+\be-\sqrt{\De})^2-v^2}{4\De}.
\end{align*}
\end{lemma}

\subsection{Cyclotomic Classes, Additive Characters, and Multiplicative Characters over Finite Fields}

Throughout the rest of this paper, we write $q=p^s$ for some prime $p$ and $s \ge 1$. Let $\Fq$ be the finite field of order $q$ and $\Fq^*$ be the multiplicative subgroup of $\Fq$.

Let $e>1$ be an integer such that $e \mid (q-1)$. Let $\om$ be a primitive element of $\Fq$. Then $\Fq$ contains a unique multiplicative subgroup of order $\frac{q-1}{e}$, written as $C_0^{(e,q)}=\pointy{\om^e}$. Then for $0 \le i \le e-1$, the multiplicative cosets $C_{i}^{(e,q)}=\om^i\pointy{\om^e}$ are the \emph{cyclotomic classes} of order $e$ over $\Fq$. We write $C_0=C_0^{(e,q)}$ and define 
$$C_0^{\perp}=\{ \chi \in \wh{\Fq^*} \mid  \mbox{ $\chi$ is principal on $C_0$}\}.
$$ 
For a standard reference on cyclotomic classes, please refer to \cite{Storer}.

Since finite fields have additive and multiplicative group structure, we shall describe characters over both the additive and the multiplicative group of finite fields.  For $n \ge 1$, we use $\xi_n$ to denote a primitive $n$-th root of unity. Below, we will define and describe some basic properties of trace functions, additive characters, and multiplicative characters over finite field. For a comprehensive treatment on these objects, please refer to \cite[Chapter 5]{LN}.

\begin{fact}[trace and norm functions over finite fields]
Let  $m \ge 1$. The \emph{trace function} $\Trqmq$ from $\Fqm$ to $\Fq$ is defined as
$$
\Trqmq(x)=x+x^q+\cdots+x^{q^{m-1}}, \;\forall x \in \Fqm.
$$
The \emph{norm function} from $\Fqm$ to $\Fq$ is defined as
$$
\Nqmq(x)=x \cdot x^q \cdots x^{q^{m-1}}=x^{\frac{q^m-1}{q-1}}, \forall x \in \Fqm.
$$
\end{fact}

\begin{fact}[additive characters over finite fields]
Suppose $q=p^s$ for some prime $p$ and $s \ge 1$. Then for each $a \in \Fq$, the mapping $\psi_a: x \rightarrow \xi_p^{\Trqp(ax)}$ is an \emph{additive character} over $\Fq$. The character $\psi_a$ is principal if and only if $a=0$. Moreover, the character group $\wh{\Fq}$ of the additive group $(\Fq,+)$ is precisely $\{\psi_a \mid a \in \Fq\}$. 
Recall that
$$
\de_0(a)=\begin{cases}
                   1 & \mbox{if $a=0$,} \\
                   0 & \mbox{if $a \ne 0$.}
                   \end{cases}
$$
When applying an additive character $\psi_a$ to the group $(\Fq,+)$, we have the following \emph{orthogonal relation}:
\begin{equation}
\label{eqn-addorthogonal}
\psi_a(\Fq)=q\de_0(a)=\begin{cases}
                   q & \mbox{if $a=0$,} \\
                   0 & \mbox{if $a \ne 0$.}
                   \end{cases}
\end{equation}  
\end{fact} 

\begin{fact}[multiplicative characters over finite fields]
 Let $\om$ be a primitive element of $\Fq$. Then every nonzero element of $\Fq^*$ can be written as $\om^i$ for some $0 \le i \le q-2$. For some $0 \le j \le q-2$, the mapping $\chi:=\chi(j)$ such that for each $0 \le i \le q-2$, $\chi(\om^i)=\xi_{q-1}^{ij}$ is a \emph{multiplicative character} over $\Fq$. The multiplicative character $\chi$ has order $\frac{q-1}{\gcd(j,q-1)}$. Let $e>1$ and $e \mid (q-1)$. Let $\chi$ be an order $e$ multiplicative character over $\Fq$. Then for an integer $j$,
\begin{equation}
\label{eqn-mulorthogonal}
\sum_{i=0}^{e-1} \chi^j(\om^i)=\begin{cases}
                   e & \mbox{if $j \equiv 0 \pmod{e}$,} \\
                   0 & \mbox{if $j \not\equiv 0 \pmod{e}$.}
                   \end{cases}
\end{equation}  
\end{fact} 

Regarding the multiplicative characters, we have the following critical notion of \emph{lifting characters} (see for instance \cite[Theorem 6]{OttDCC}).

\begin{lemma}[lifting characters]
\label{lem-lifting}
Let $m \ge 1$ and $e \mid (q-1)$. Suppose $\chi^{\pr}$ is a multiplicative character of order $e$ over $\Fqm$. Then there exists some multiplicative character $\chi$ of order $e$ over $\Fq$, such that $\chi^{\pr}=\chi \circ \Nqmq$, namely, $\chi^{\pr}(\ga)=\chi(\Nqmq(\ga))$ for each $\ga \in \Fqm^*$. In this case, $\chi^{\pr}$ is called the \emph{lifting character} of $\chi$.
\end{lemma}

\subsection{Gauss Sums}

Gauss sums were first introduced by Gauss~\cite{Gauss}, showcasing an elegant way integrating additive and multiplicative characters over finite fields. 

\begin{definition}[Gauss sum]
Let $\chi$ be a multiplicative character of $\Fq$. Let $\psi_1$ be the additive character of $\Fq$. The Gauss sum with respect to $\chi$ is defined as
$$
G(\chi)=\sum_{x \in \Fq^*} \chi(x)\psi_1(x).
$$ 
\end{definition}

Below, we record some useful properties of Gauss sums. For a comprehensive treatment of Gauss sums, please refer to \cite{BEW}.

\begin{lemma}[{\cite[Section 2]{MX}}]
\label{lem-Gauss}
Let $\chi$ be a multiplicative character over $\Fq$.  Then the following holds.
\begin{itemize}
\item[(1)] $G(\chi)=-1$ if $\chi$ is principal.
\item[(2)] $G(\chi)\ol{G(\chi)}=q$ if $\chi$ is nonprincipal.
\item[(3)] $G(\chi^{-1})=\chi(-1)\ol{G(\chi)}$.
\item[(4)] Let $\om$ be a primitive element of $\Fq$. Let $e \mid (q-1)$. Let $I$ be a set of integers that are pairwise distinct modulo $e$. Then $D=\cup_{i \in I} C_i^{(e,q)}$ is a subset of $\Fq$ being a union of cyclotomic classes of order $e$ over $\Fq$. Moreover, for $a \in \Fq^*$,
\begin{equation}
\label{eqn-gaussperiod}
\psi_a(D)=\frac{1}{e}\sum_{\chi \in C_0^{\perp}}G(\chi^{-1})\sum_{i \in I} \chi(a\om^i).
\end{equation}
\end{itemize}
\end{lemma}

The following Davenport-Hasse lifting formula deals with Gauss sums with respect to the lifting characters.

\begin{lemma}[Davenport-Hasse]
\label{lem-DH}
Let $\chi$ be a nonprincipal multiplicative character of $\Fq$. Let $m \ge 1$ and $\chi^{\pr}$ be the lifting character of $\chi$ over $\Fqm$, namely, $\chi^{\pr}(\ga)=\chi(\Nqmq(\ga))$ for each $\ga \in \Fqm^*$. Then
$$
G_m(\chi^{\pr})=(-1)^{m-1}G_1(\chi)^m,
$$
where $G_m(\chi^{\pr})$ is the Gauss sum with respect to $\chi^{\pr}$ over $\Fqm$ and $G_1(\chi)$ is the Gauss sum with respect to $\chi$ over $\Fq$.
\end{lemma}

\section{Statement of the Main Result}
\label{sec-statement}

In this section, we present explicit constructions of SRGs and PDSs that give rise to Theorems~\ref{theorem-main} and \ref{theorem-dualmain}. Moreover, we examine how our constructions relate to previously known results. Our construction relies on cyclotomic classes over finite fields, which have been an extremely powerful tool in constructing SRGs and PDSs, see for instance \cite{BXZ,BWX,DHJP,M,MX,OttDM,P}. 

We first give an explicit construction of the family of SRGs and PDSs in Theorem~\ref{theorem-main}.

\begin{theorem}
\label{theorem-main2}
Let $s \ge 1$ and $q=p^s$ be a power of prime $p$. Let $m$ and $\ell$ be positive integers. Set $e=\frac{q^m-1}{q-1}$ and $\ga$ to be a primitive element of $\Fqm$. For $0 \leq r \leq m$, let $R$ be an $r$-dimensional vector subspace of $\Fqm$ over $\Fq$. Define 
\begin{equation}
\label{eqn-T}
T=\{ 0 \le i < e \mid \ga^i \in R \}.
\end{equation}
Then 
$$
D_{r,T}=\bigcup_{i=0}^{e-1} (\Cqml_i \times \cup_{t \in T} \Cqmlo_{i+t}) \bigcup (\Fqml^* \times \{ 0 \})
$$ 
is a $(v,k_r,\la_r,\mu_r)$-PDS in the elementary abelian $p$-group $G=(\Fqml \times \Fqmlo,+) \cong \Z_p^{sm(2\ell+1)}$, where
\begin{align*}
v = & q^{m(2 \ell+1)}, \\
k_r  = & \frac{(q^r-1)(q^{m\ell}-1)(q^{m(\ell+1)}-1)}{q^m-1} + q^{m\ell}-1, \\
\lambda_r = & \frac{(q^r-1)(q^{m(\ell+1)}-1)}{q^m-1}\Big(\frac{(q^r-1)(q^{m\ell}-1)}{q^m-1}-1\Big) + q^{m\ell}-2, \\
\mu_r = & \frac{(q^r-1)(q^{m\ell}-1)}{q^m-1}\Big(\frac{(q^r-1)(q^{m(\ell+1)}-1)}{q^m-1}+1\Big). 
\end{align*}
Alternatively, we have
\begin{equation}
\label{eqn-DTalt}
\begin{split}
D_{r,T} =&\{ (a,b) \mid ab \ne 0, \frac{\Nqmloqm(b)}{\Nqmlqm(a)} \in \cup_{t \in T} \Cqm_t \} \cup \{ (a,0) \mid a \in \Fqml^*\} \\ 
            =&\{ (a,b) \mid ab \ne 0, \frac{\Nqmloqm(b)}{\Nqmlqm(a)} \in R^* \} \cup \{ (a,0) \mid a \in \Fqml^*\}  
\end{split}            
\end{equation}
Consequently, $\Cay(G,D_{r,T})$ is a $(v,k_r,\la_r,\mu_r)$-SRG with a maximum clique $\Fqml \times \{ 0 \}$ of size $q^{m\ell}$.
\end{theorem}

A few remarks regarding the family of SRGs and PDSs derived from Theorem~\ref{theorem-main2} are in order.
\begin{remark}
\label{remark-main2}
\begin{itemize}
\item[(1)] The cases of $r \in \{0,m \}$ in Theorem~\ref{theorem-main2} correspond to degenerate SRGs and PDSs. Indeed, when $r=0$, we have $T=\es$. Then $D_{0,\es}=\Fqml^* \times \{0\}$ and $\Cay(G,D_{0,\es})$ is a union of $q^{m(\ell+1)}$ copies of the complete graph $K_{q^{m\ell}}$. When $r=m$, we have $T=\{0,1,\ldots,e-1\}$. Then $D_{m,\{0,1,\ldots,e-1\}}=\Fqml^* \times \Fqmlo$ and $\Cay(G,D_{m,\{0,1,\ldots,e-1\}})$ is a complete $q^{m\ell}$-partite graph 
$$
K_{ \underbrace{q^{m(\ell+1)},\ldots,q^{m(\ell+1)}}_\text{$q^{m\ell}$} }
$$ 
with each part having size $q^{m(\ell+1)}$.
\item[(2)] In \cite[Theorem 3]{OttDM}, Ott constructed a family of PDSs in $G=(\Fqml \times \Fqmlo,+)$. Specifically, let $\chi$ be a multiplicative character of order $\frac{q^m-1}{q-1}$  over $\Fqm$. Let $\psi$ be the lifting character of $\chi^{-1}$ over $\Fqml$ and $\vp$ be the lifting character of $\chi$ over $\Fqmlo$. Then 
$$
D=\{ (a,b) \mid ab \ne 0,  \psi(a)\vp(b)=1\} \cup \{ (a,0) \mid a \in \Fqml^*\}
$$
is a $(v,k_1,\la_1,\mu_1)$-PDS in $G$. Indeed, we can rephrase $D$ as follows:
$$
D=\{ (a,b) \mid ab \ne 0,  \frac{\Nqmloqm(b)}{\Nqmlqm(a)} \in \Cqm_0\} \cup \{ (a,0) \mid a \in \Fqml^*\}.
$$
Therefore, applying Theorem~\ref{theorem-main2} with $R$ being an $1$-dimensional subspace of $\Fqm$ over $\Fq$ and $T=\{0\}$, we recover \cite[Theorem 3]{OttDM}.
\item[(3)] In \cite{OttDCC}, Ott presented a very broad framework in order to construct PDSs in the additive group of Galois domains, which are direct product of finite fields. Specifically, Ott gave the following construction.
\begin{res}[{\cite[Theorem 20]{OttDCC}}]
\label{res-Ott}
Let $q_1$, $q_2$ be prime powers. Let $d$ be a positive integer such that $d \mid \gcd(q_1-1,q_2-1)$ and $1 \le u \le d-1$. Let $\De=\Ga \cup (\F_{q_1}^* \times \{ 0 \})$. Then $\De$ is a $(v,k,\la.\mu)$-PDS in $(\F_{q_1} \times \F_{q_2},+)$ with 
\begin{align*}
v =& q_1q_2, \\
k =& \frac{u(q_1-1)(q_2-1)}{d}+q_1-1, \\
\la =& \frac{u(q_2-1)}{d}(\frac{u(q_1-1)}{d}-1)+q_1-2, \\
\mu=& \frac{u(q_1-1)}{d}(\frac{u(q_2-1)}{d}+1),
\end{align*}
if $\Ga=\{ (-a,-b) \mid (a,b) \in \Ga \}$ and the following holds true:
\begin{itemize}
\item[(a)] $\Ga$ is a $(v,k-q_1+1,\la^*,\mu^*)$-PDS in $(\F_{q_1} \times \F_{q_2},+)$,
\item[(b)] $\la^*-\mu^*=\la-\mu+2=q_1-\frac{u(q_1+q_2-2)}{d}$.
\end{itemize}
\end{res}
On one hand, set $q_1=\qml$, $q_2=\qmlo$, $d=\frac{q^m-1}{q-1}$ and $u=\frac{q^r-1}{q-1}$, assuming all the conditions in Result~\ref{res-Ott} hold true, then the corresponding PDS $\De$ has the same parameters as $D_{r,T}$ in Theorem~\ref{theorem-main2}. On the other hand, the construction in Theorem~\ref{theorem-main2} does not belong to the framework in Result~\ref{res-Ott}. Assuming otherwise, the subset $\bigcup_{i=0}^{e-1} (\Cqml_i \times \cup_{t \in T} \Cqmlo_{i+t})$ of $D_{r,T}$ would correspond to the PDS $\Ga$ in Result~\ref{res-Ott}. However, by Equation~(\ref{eqn-DT4}), 
\begin{align*} 
    & (\psi_a \times \vp_b)(\bigcup_{i=0}^{e-1} (\Cqml_i \times \cup_{t \in T} \Cqmlo_{i+t})) \\
 = & \begin{cases} 
          -\frac{(q^r-1)(q^{m(\ell+1)}-1)}{q^m-1} & \mbox{if $a \ne 0$, $b=0$ or $ab \ne 0$, $\frac{\Nqmloqm(b)}{\Nqmlqm(a)} \in R^{\perp}$,} \\
           -\frac{(q^r-1)(q^{m\ell}-1)}{q^m-1} & \mbox{if $a=0$, $b \ne 0$,} \\
           \frac{(q^m-q^r)(q^{m\ell}-1)}{q^m-1}+1 & \mbox{if $ab \ne 0$, $\frac{\Nqmloqm(b)}{\Nqmlqm(a)} \not\in R^{\perp}$.}
        \end{cases}       
\end{align*}
Consequently, $\bigcup_{i=0}^{e-1} (\Cqml_i \times \cup_{t \in T} \Cqmlo_{i+t})$ is not a PDS, leading to a contradiction.
\item[(4)] In \cite[Theorem 4.1]{DeWinter}, De Winter described an algebraic construction of a projective set, which can be rephrased into the context of PDS as follows. Suppose that $\epsilon$ is a primitive element of $\F_{q^{2m}}$, then $\tau=\epsilon^{q^m+1}$ is a primitive element of $\Fqm$. Then the subset 
\begin{align*} 
D=&\{ (\tau^i,\epsilon^i\epsilon^{j\frac{q^m-1}{q-1}}) \mid 0 \le i < \frac{q^m-1}{q-1}, 0 \le j < (q-1)(q^m+1)-1 \} \\
                                                                                                                                    & \cup (\Fqm^* \times \{ 0 \}) \\
   =&\{ (a,b) \mid ab \ne 0,  \frac{N_{q^{2m}/q^m}(b)}{a} \in C^{(\frac{q^m-1}{q-1},q^m)}_0 \} \cup \{ (a,0) \mid a \in \Fqm^*\}
\end{align*}
is a $(q^{3m},(q^{m+1}-q^m+q)(q^m-1),q^m-q+(q^{m+1}-q^m+q)(q-2),(q^{m+1}-q^m+q)(q-1))$-PDS in $\Fqm \times \F_{q^{2m}}$. In view of Equation~(\ref{eqn-DTalt}), Theorem~\ref{theorem-main2} indeed generalizes \cite[Theorem 4.1]{DeWinter}. Moreover, an open problem in \cite[Section 5]{DeWinter} essentially asks if it is possible to glue multiple PDSs following from \cite[Theorem 4.1]{DeWinter}, where they all share a common subset $\Fqm^* \times \{ 0 \}$, so that the union of these PDSs is again a PDS in $\Fqm \times \F_{q^{2m}}$. Theorem~\ref{theorem-main2} provides an affirmative answer to this question.  
\item[(5)] In view of Equation~(\ref{eqn-DTalt}), for the primitive element $\ga$ of $\Fqm$, we have
\begin{align*}
\ga D_{r,T} =&\{ (\ga a,\ga b) \mid ab \ne 0, \frac{\Nqmloqm(b)}{\Nqmlqm(a)} \in R^* \} \cup \{ (\ga a,0) \mid a \in \Fqml^*\} \\ 
            =&\{ (a,b) \mid ab \ne 0, \frac{\Nqmloqm(\ga^{-1}b)}{\Nqmlqm(\ga^{-1}a)} \in R^* \} \cup \{ (a,0) \mid a \in \Fqml^*\}  \\
            =&\{ (a,b) \mid ab \ne 0, \ga^{-1}\frac{\Nqmloqm(b)}{\Nqmlqm(a)} \in R^* \} \cup \{ (a,0) \mid a \in \Fqml^*\} \\
            =&\{ (a,b) \mid ab \ne 0, \frac{\Nqmloqm(b)}{\Nqmlqm(a)} \in \ga R^* \} \cup \{ (a,0) \mid a \in \Fqml^*\} 
\end{align*}            
If $r=0$, then $\ga R^*=R^*=\es$ and $T=\es$, which implies $\ga D_{0,\es}=D_{0,\es}=\Fqml^* \times \{0\}$. If $r=m$, then $\ga R^*=R^*=\Fqm^*$ and $T=\{0,1,\ldots,e-1\}$, which implies $\ga D_{m,\{0,1,\ldots,e-1\}}=D_{m,\{0,1,\ldots,e-1\}}=\Fqml^* \times \Fqmlo$. If $1 \le r \le m-1$, then $\ga R^* \ne R^*$ and therefore, $\ga D_{r,T} \ne D_{r,T}$. Thus, $D_{r,T}$ is not invariant under the multiplication of elements from $\Fqm^*$. Consequently, the projective set in $\PG(m(2\ell+1)-1,q)$ corresponding to $D_{r,T}$ cannot be the blow up of a point set in $\PG(2,q^m)$.
\end{itemize}
\end{remark}

Next, we describe an explicit construction of the family of SRGs and PDSs with dual generalized Denniston parameters in Theorem~\ref{theorem-dualmain}.

\begin{theorem}
\label{theorem-dualmain2}
Let $s \ge 1$ and $q=p^s$ be a power of prime $p$. Let $m$ and $\ell$ be positive integers. Set $e=\frac{q^m-1}{q-1}$ and $\ga$ to be a primitive element of $\Fqm$. For $0 \leq r \leq m$, let $R$ be an $r$-dimensional vector subspace of $\Fqm$ over $\Fq$, which is also an $sr$-dimensional vector subspace of $\Fqm$ over $\Fp$. Let 
$$
R^{\perp} = \{ x \in \Fqm \mid \Trqmp(xy) = 0 \mbox{ for each } y \in R \}
$$
be the dual space of $R$ with respect to the trace function $\Trqmp$, which is an $s(m-r)$-dimensional vector subspace space of $\Fqm$ over $\Fp$. Define 
\begin{equation*}
T^{\perp}=\{ 0 \le i < e \mid \ga^i \in R^{\perp} \}
\end{equation*}
and
$$
{T^{\perp}}^c=\{ 0,1,\ldots,e-1 \} \sm T^{\perp}.
$$
Then 
\begin{equation}
\label{eqn-DTp}
D_{r,T}^+=\bigcup_{i=0}^{e-1} (\Cqml_i \times \cup_{t \in {T^{\perp}}^c} \Cqmlo_{i+t}) \bigcup (\{0\} \times \Fqmlo^*)
\end{equation}
is a $(v^+,k_r^+,\la_r^+,\mu_r^+)$-PDS in the elementary abelian $p$-group $G=(\Fqml \times \Fqmlo,+) \cong \Z_p^{sm(2\ell+1)}$, where
\begin{equation*}
\begin{split}
v^+ = & q^{m(2 \ell+1)}, \\
k_r^+  = & \frac{(q^m-q^{m-r})(q^{m\ell}-1)(q^{m(\ell+1)}-1)}{q^m-1} + q^{m(\ell+1)}-1, \\
\lambda_r^+ = & \frac{(q^m-q^{m-r})(q^{m(\ell+1)}-1)(q^{m(\ell+1)}-q^{m(\ell+1)-r}+q^{m-r}-1)}{(q^m-1)^2} \\
                        & + \frac{q^{m(\ell+2)-r}+q^{m(\ell+1)-r}-2q^{m(\ell+1)}-2q^{m-r}+2}{q^m-1}, \\
\mu_r^+ = & \frac{(q^m-q^{m-r})(q^{m(\ell+1)}-1)(q^{m(\ell+1)}-q^{m(\ell+1)-r}+q^{m-r}-1)}{(q^m-1)^2}.
\end{split}
\end{equation*}
Alternatively, we have
\begin{equation}
\label{eqn-DTpalt}
\begin{split}
D_{r,T}^+ =&\{ (a,b) \mid ab \ne 0, \frac{\Nqmloqm(b)}{\Nqmlqm(a)} \in \cup_{t \in {T^\perp}^c} \Cqm_t \} \cup \{ (0,b) \mid b \in \Fqmlo^*\} \\
            =&\{ (a,b) \mid ab \ne 0, \frac{\Nqmloqm(b)}{\Nqmlqm(a)} \notin R^\perp \} \cup \{ (0,b) \mid b \in \Fqmlo^*\} 
\end{split}            
\end{equation}
Consequently, $\Cay(G,D_{r,T}^+)$ is a $(v^+,k_r^+,\la_r^+,\mu_r^+)$-SRG with a maximum clique $\{0\} \times \Fqmlo$ of size $q^{m(\ell+1)}$.
\end{theorem}

\begin{remark}
From Equation~(\ref{eqn-DTp}), we can see that $D_{r,T}^+$ is the complement of $D_{m-r,T^{\perp}}$, where the latter has the generalized Denniston paramters $(v,k_{m-r},\la_{m-r},\mu_{m-r})$. 
\end{remark}

\section{Proof of the Main Result}
\label{sec-proof}

Now we proceed to prove Theorem~\ref{theorem-main2}, which implies Theorem~\ref{theorem-main}.

\begin{proof}[Proof of Theorem~\ref{theorem-main2}]
Note that $|T| = \frac{q^r-1}{q-1}$. Then 
$$
k_r=|D_{r,T}|=\frac{(q^r-1)(q^{m\ell}-1)(q^{m(\ell+1)}-1)}{q^m-1} + q^{m\ell}-1. 
$$
Note that 
$k_r^2=k_r-\mu_r+(\la_r-\mu_r)k_r+\mu_r v$ and $(\la_r-\mu_r)^2+4(k_r-\mu_r) \ge 0$. In view of Lemma~\ref{lem-pdschar}, it suffices to show that
$$
\{ \chi(D_{r,T}) \mid \mbox{$\chi \in \wh{G}$ nonprincipal} \}=\{ -\frac{(q^r-1)(q^{m(\ell+1)}-1)}{q^m-1}-1, \frac{(q^m-q^r)(q^{m\ell}-1)}{q^m-1} \}.
$$

Note that $G=(\Fqml \times \Fqmlo,+)$, then $\wh{G}=\{ \psi_a \times \vp_b \mid a \in \Fqml, b \in \Fqmlo \}$, where
\begin{align*}
\psi_a(x_1) = & \xi_p^{\Trqmlp(ax_1)}, \forall x_1 \in \Fqml, \\
\vp_b(x_2)  = & \xi_p^{\Trqmlop(bx_2)}, \forall x_2 \in \Fqmlo, \\
(\psi_a \times \vp_b)(x_1,x_2)  = & \psi_a(x_1)\vp_b(x_2), \forall (x_1, x_2) \in G.
\end{align*}
Together with Equation~(\ref{eqn-addorthogonal}), we have
\begin{align}
(\psi_a \times \vp_b)(D_{r,T})=&\sum_{i=0}^{e-1} \psi_a(\Cqml_i) \vp_b(\cup_{t \in T} \Cqmlo_{i+t})+\psi_a(\Fqml^*) \nonumber \\
                                          =&\sum_{i=0}^{e-1} \psi_a(\Cqml_i) \sum_{t \in T} \vp_b(\Cqmlo_{i+t})+q^{m\ell}\de_0(a)-1
\label{eqn-DT1}
\end{align}

We shall consider three cases.

{\flushleft {\bf Case 1}: $a=0$ and $b \neq 0$, namely, $\psi_a$ is principal and $\vp_b$ is nonprincipal.}

In this case, $\psi_0(\Cqml_i) = |\Cqml_i| = \frac{(q-1)(q^{m\ell}-1)}{q^m-1}$ for each $0 \le i \le e-1$. Note that for each $t \in T$, $\cup_{i=0}^{e-1} \Cqmlo_{i+t}=\Fqmlo^*$. By Equation~(\ref{eqn-DT1}) and Equation~(\ref{eqn-addorthogonal}),
\begin{align*}
(\psi_0 \times \vp_b)(D_{r,T})  = & \frac{(q-1)(q^{m\ell}-1)}{q^m-1} \sum_{t \in T} \sum_{i=0}^{e-1} \vp_b(\Cqmlo_{i+t})+ q^{m\ell}-1 \\
 = & \frac{(q-1)(q^{m\ell}-1)}{q^m-1} \sum_{t \in T} \vp_b(\cup_{i=0}^{e-1} \Cqmlo_{i+t})+ q^{m\ell}-1 \\
 = & \frac{(q-1)(q^{m\ell}-1)}{q^m-1} \sum_{t \in T} \vp_b(\Fqmlo^*) + q^{m\ell}-1 \\
 = & \frac{(q-1)(q^{m\ell}-1)}{q^m-1} |T| \vp_b(\Fqmlo^*) + q^{m\ell}-1 \\
 = & -\frac{(q-1)(q^{m\ell}-1)}{q^m-1} \cdot \frac{q^r-1}{q-1} + q^{m\ell}-1 \\
 = & \frac{(q^m-q^r)(q^{m\ell}-1)}{q^m-1}
\end{align*}

{\flushleft {\bf Case 2}: $a \ne 0$ and $b=0$, namely, $\psi_a$ is nonprincipal and $\vp_b$ is principal.}

In this case, $\vp_0(\Cqmlo_i) =  |\Cqmlo_i| = \frac{(q-1)(q^{m(\ell+1)}-1)}{q^{m}-1}$ for each $0 \le i \le e-1$. Note that $\cup_{i=0}^{e-1} \Cqml_{i}=\Fqml^*$. By Equation~(\ref{eqn-DT1}),
\begin{align*}
(\psi_a \times \vp_0)(D_{r,T}) = & \sum_{i=0}^{e-1} \psi_a(\Cqml_i) |T| \frac{(q-1)(q^{m(\ell+1)}-1)}{q^{m}-1} - 1 \\
 = & \frac{q^r-1}{q-1} \cdot \frac{(q-1)(q^{m(\ell+1)}-1)}{q^{m}-1} \psi_a(\cup_{i=0}^{e-1} \Cqml_i) - 1 \\
 = & \frac{(q^r-1)(q^{m(\ell+1)}-1)}{q^{m}-1} \psi_a(\Fqml^*) - 1 \\
 = & -\frac{(q^r-1)(q^{m(\ell+1)}-1)}{q^m-1}-1
\end{align*}

{\flushleft {\bf Case 3}: $ab \ne 0$, namely, $\psi_a$ and $\vp_b$ are both nonprincipal.}

Let $\al$ be a primitive element of $\Fqml$. For $0 \le i \le e-1$, set $\Cqml_i=\al^i\pointy{\al^e}$ and we write $a \in \Fqml^*$ as $a = \al^{f_a}$ for some $0 \le f_a \le q^{m\ell}-2$. Note that $\Nqmlqm(\al)=\al^{\frac{\qml-1}{q^m-1}}$ is a primitive element of $\Fqm$, we can choose $\al$ such that $\Nqmlqm(\al)=\ga$, where $\ga$ is the primitive element of $\Fqm$ used in Equation~(\ref{eqn-T}) to define the subset $T$.

Similarly,  we can choose a primitive element $\be$ of $\Fqmlo$, such that  $\Nqmloqm(\be)=\ga$. For $0 \le i \le e-1$, set $\Cqmlo_i=\be^i\pointy{\be^e}$ and we write $b \in \Fqmlo^*$ as $b = \be^{f_b}$ for some $0 \le f_b \le q^{m(\ell+1)}-2$. By Equation~(\ref{eqn-DT1}),
\begin{align}
\label{eqn-DT2}
(\psi_a \times \vp_b)(D_{r,T})+1= & \sum_{i=0}^{e-1} \psi_{\al^{f_a}}(\Cqml_i) \sum_{t \in T} \vp_{\be^{f_b}}(\Cqmlo_{i+t}) \nonumber \\
                                              = & \sum_{i=0}^{e-1} \psi_{1}(\Cqml_{i+f_a}) \sum_{t \in T} \vp_{1}(\Cqmlo_{i+f_b+t})
\end{align}
For simplicity, we write
\begin{align*}
\Cl = & {\Cqml_0}^{\perp}=\{ \chi^{\pr} \in \widehat{\Fqml^*} \mid \mbox{$\chi^{\pr}$ is principal on $\Cqml_0$} \}, \\
\Clo = & {\Cqmlo_0}^{\perp}= \{ \tht^{\pr} \in \widehat{\Fqmlo^*} \mid \mbox{$\tht^{\pr}$ is principal on $\Cqmlo_0$} \}.
\end{align*}
Note that $\Cl \cong \Fqml^*/\Cqml_0$ and $\Clo \cong \Fqmlo^*/\Cqmlo_0$. Both $\Cl$ and $\Clo$ are isomorphic to the cyclic group of order $e$. Assume that $\chi^{\pr}_e$ and $\tht^{\pr}_e$ are order $e$ characters in $\Cl$ and $\Clo$ respectively. Then
\begin{equation}
\label{eqn-chargp}
\Cl = \pointy{\chi^{\pr}_e}=\{ {\chi^{\pr}_e}^y \mid 0 \le y \le e-1 \} \mbox{ and }\Clo = \pointy{\tht^{\pr}_e}= \{ {\tht^{\pr}_e}^z \mid 0 \le z \le e-1 \}.
\end{equation}
According to Lemma~\ref{lem-lifting}, since $e=\frac{q^m-1}{q-1} \mid (q^m-1)$, both $\chi^{\pr}_e$ and $\tht^{\pr}_e$ are lifting characters of certain order $e$ multiplicative characters over $\Fqm$. Note that all order $e$ multiplicative characters over $\Fqml$ are of the form ${\chi^{\pr}_e}^y$ with $0 \le y \le e-1$ and $\gcd(y,e)=1$, as well as all order $e$ multiplicative characters over $\Fqmlo$ are of the form ${\tht^{\pr}_e}^z$ with $0 \le z \le e-1$ and $\gcd(z,e)=1$. Therefore, we can choose order $e$ characters $\chi^{\pr}_e$ and $\tht^{\pr}_e$ such that they are the lifting characters of the same order $e$ multiplicative character $\chi_e$ over $\Fqm$. Namely, 
$$
\chi^{\pr}_e = \chi_e \circ \Nqmlqm  \mbox{ and } \tht^{\pr}_e = \chi_e \circ \Nqmloqm.
$$
Consequently, we have
\begin{equation}
\label{eqn-lifting}
\chi^{\pr}_e(\al)=\chi_e \circ \Nqmlqm(\al)=\chi_e(\ga) \mbox{ and } \tht^{\pr}_e(\be) = \chi_e \circ \Nqmloqm(\be)=\chi_e(\ga).
\end{equation}
Let $G_{m\ell}(\chi^{\pr})$ be the Gauss sum with respect to the multiplicative character $\chi^{\pr}$ over $\Fqml$. Let $G_{m(\ell+1)}(\tht^{\pr})$ be the Gauss sum with respect to the multiplicative character $\tht^{\pr}$ over $\Fqmlo$.  Let $G_{m}(\chi)$ be the Gauss sum with respect to the multiplicative character $\chi$ over $\Fqm$.
By Equations~(\ref{eqn-gaussperiod}), (\ref{eqn-chargp}), (\ref{eqn-lifting}), Lemma~\ref{lem-Gauss}(1) and the Davenport-Hasse formula, we have 
\begin{align*}
\psi_1(\Cqml_{i+f_a})= & \frac{1}{e} \sum_{\chi^{\pr} \in \Cl} G_{m\ell}({\chi^{\pr}}^{-1})\chi^{\pr}(\al^{i+f_a}) \\
                                  = & \frac{1}{e} \sum_{y=0}^{e-1} G_{m\ell}({\chi_e^{\pr}}^{-y}){\chi_e^{\pr}}^y(\al^{i+f_a}) \\
                                  = & \frac{1}{e} \Big(\sum_{y=1}^{e-1} G_{m\ell}({\chi_e^{\pr}}^{-y}){\chi_e^{\pr}}^y(\al^{i+f_a})-1 \Big) \\
                                  = & \frac{1}{e} \Big(\sum_{y=1}^{e-1} (-1)^{\ell-1} G_{m}(\chi_e^{-y})^{\ell} \chi_e^y(\ga^{i+f_a})-1 \Big)
\end{align*}
and
\begin{align*}                                  
\vp_1(\Cqmlo_{i+f_b+t}) = & \frac{1}{e} \sum_{\tht^{\pr} \in \Clo} G_{m(\ell+1)}({\tht^{\pr}}^{-1}) \tht^{\pr}(\be^{i+f_b+t}) \\ 
                                       = & \frac{1}{e} \sum_{z=0}^{e-1} G_{m(\ell+1)}({\tht_e^{\pr}}^{-z}) {\tht_e^{\pr}}^z(\be^{i+f_b+t}) \\
                                       = & \frac{1}{e} \Big(\sum_{z=1}^{e-1} G_{m(\ell+1)}({\tht_e^{\pr}}^{-z}) {\tht_e^{\pr}}^z(\be^{i+f_b+t})-1 \Big) \\
                                       = & \frac{1}{e} \Big(\sum_{z=1}^{e-1} (-1)^{\ell}G_{m}(\chi_e^{-z})^{\ell+1} \chi_e^z(\ga^{i+f_b+t})-1 \Big)
\end{align*}
By the above two equations and Equation~(\ref{eqn-DT2}), we have
\begin{align*}
   & e^2((\psi_a \times \vp_b)(D_{r,T})+1)  \\
= & \sum_{i=0}^{e-1} \Big(\sum_{y=1}^{e-1} (-1)^{\ell-1} G_{m}(\chi_e^{-y})^{\ell} \chi_e^y(\ga^{i+f_a})-1 \Big) \sum_{t \in T} \Big(\sum_{z=1}^{e-1} (-1)^{\ell}G_{m}(\chi_e^{-z})^{\ell+1} \chi_e^z(\ga^{i+f_b+t})-1 \Big) \\
= & \sum_{i=0}^{e-1} \Big(\sum_{y=1}^{e-1} (-1)^{\ell-1} G_{m}(\chi_e^{-y})^{\ell} \chi_e^y(\ga^{i+f_a})-1 \Big)  \Big( \sum_{t \in T} \sum_{z=1}^{e-1} (-1)^{\ell}G_{m}(\chi_e^{-z})^{\ell+1} \chi_e^z(\ga^{i+f_b+t})-|T| \Big) \\
= & \sum_{i=0}^{e-1} \Big(-\sum_{y=1}^{e-1} G_{m}(\chi_e^{-y})^{\ell} \chi_e^y(\ga^{i+f_a}) \sum_{z=1}^{e-1} \sum_{t \in T} G_{m}(\chi_e^{-z})^{\ell+1} \chi_e^z(\ga^{i+f_b+t}) \\
   & - |T| \sum_{y=1}^{e-1} (-1)^{\ell-1} G_{m}(\chi_e^{-y})^{\ell} \chi_e^y(\ga^{i+f_a}) - \sum_{t \in T} \sum_{z=1}^{e-1} (-1)^{\ell}G_{m}(\chi_e^{-z})^{\ell+1} \chi_e^z(\ga^{i+f_b+t}) +|T| \Big) \\
= & - \sum_{y=1}^{e-1} G_{m}(\chi_e^{-y})^{\ell} \chi_e^y(\ga^{f_a}) \sum_{z=1}^{e-1} \sum_{t \in T} G_{m}(\chi_e^{-z})^{\ell+1} \chi_e^z(\ga^{f_b+t})\sum_{i=0}^{e-1}  \chi_e^{y+z}(\ga^{i}) \\
   & - |T| \sum_{y=1}^{e-1} (-1)^{\ell-1} G_{m}(\chi_e^{-y})^{\ell} \chi_e^y(\ga^{f_a}) \sum_{i=0}^{e-1} \chi_e^y(\ga^{i}) \\
   & - \sum_{z=1}^{e-1} \sum_{t \in T} (-1)^{\ell}G_{m}(\chi_e^{-z})^{\ell+1} \chi_e^z(\ga^{f_b+t}) \sum_{i=0}^{e-1} \chi_e^z(\ga^{i}) +e|T| 
\end{align*}
By Equation~(\ref{eqn-mulorthogonal}), we have
\begin{align*}
   & e^2((\psi_a \times \vp_b)(D_{r,T})+1)  \\
 =& - e \sum_{y=1}^{e-1} G_{m}(\chi_e^{-y})^{\ell} \chi_e^y(\ga^{f_a}) \sum_{t \in T} G_{m}(\chi_e^{y-e})^{\ell+1} \chi_e^{e-y}(\ga^{f_b+t})+e|T| \\
 =& - e \sum_{y=1}^{e-1} G_{m}(\chi_e^{-y})^{\ell} G_{m}(\chi_e^{y})^{\ell+1} \sum_{t \in T}  \chi_e^y(\ga^{f_a-f_b-t})+e|T| \\
\end{align*}
By Lemma~\ref{lem-Gauss}(2)(3) and Equation~(\ref{eqn-gaussperiod}), 
\begin{align*}
   & e((\psi_a \times \vp_b)(D_{r,T})+1)  \\
 =& - \sum_{y=1}^{e-1}  \chi_e^{y}((-1)^{\ell})  \ol{G_{m}(\chi_e^y)}^{\ell} G_{m}(\chi_e^{y})^{\ell+1} \sum_{t \in T}  \chi_e^y(\ga^{f_a-f_b-t})+|T| \\
 =& - \qml \sum_{y=1}^{e-1}  \chi_e^{y}((-1)^{\ell})  G_{m}(\chi_e^{y}) \sum_{t \in T}  \chi_e^y(\ga^{f_a-f_b-t})+|T| \\
 =& - \qml \sum_{y=1}^{e-1}  G_{m}(\chi_e^{y}) \sum_{t \in T}  \chi_e^{-y}((-1)^{-\ell}\ga^{t+f_b-f_a})+|T|  \\
 =& - \qml \Big(\sum_{y=0}^{e-1}  G_{m}(\chi_e^{y}) \sum_{t \in T}  \chi_e^{-y}((-1)^{\ell}\ga^{t+f_b-f_a}) + |T| \Big)+|T|  \\
 =& - \qml \sum_{y=0}^{e-1}  G_{m}(\chi_e^{y}) \sum_{t \in T}  \chi_e^{-y}((-1)^{\ell}\ga^{f_b-f_a}\ga^t) -(\qml-1)|T|
\end{align*}
Note that $\chi_e \in \wh{\Fqm^*}$ has order $e$ and  ${\Cqm_0}^{\perp}=\{ \chi_e^y \mid 0 \le y \le e-1 \}$. For $d \in \Fqm$, let $\phi_d$ be the additive character over $\Fqm$ and write $\tl=(-1)^{\ell}\ga^{f_b-f_a}$. By Equation~(\ref{eqn-gaussperiod}), we have
\begin{align*}
  (\psi_a \times \vp_b)(D_{r,T})+1 =& - \frac{\qml}{e}\sum_{\chi \in {\Cqm_0}^{\perp}}  G_{m}(\chi^{-1}) \sum_{t \in T}  \chi (\tl\ga^t) - \frac{\qml-1}{e}|T|  \\
  =& - \qml \sum_{t \in T}  \phi_{\tl} (\Cqm_t) - \frac{(q^r-1)(\qml-1)}{q^m-1}  \\
  =& - \qml \sum_{t \in T}  \phi_{1} (\tl\ga^t\Cqm_0) - \frac{(q^r-1)(\qml-1)}{q^m-1} \\
  =& - \qml \sum_{t \in T}  \phi_{1} ((-1)^{\ell}\ga^{t+f_b-f_a}\Cqm_0) - \frac{(q^r-1)(\qml-1)}{q^m-1}
\end{align*}
Crucially, since $e=\frac{q^m-1}{q-1}$, we have $\Cqm_0=\Fq^*$. Note that $(-1)^\ell \in \Fq^*$. We have
\begin{align}
  (\psi_a \times \vp_b)(D_{r,T}) =& - \qml \sum_{t \in T}  \phi_{1} (\ga^{t+f_b-f_a}\Fq^*) - \frac{(q^r-1)(\qml-1)}{q^m-1}-1  \nonumber \\
                                        =& - \qml \sum_{t \in T}  \sum_{\om\in\Fq^*}\zp^{\Trqmp(\ga^{f_b-f_a} \cdot \ga^t \om)} - \frac{(q^r-1)(\qml-1)}{q^m-1}-1       
\label{eqn-DT3}                                                                         
\end{align}
Recall that $R$ is an $r$-dimensional vector subspace of $\Fqm$ over $\Fq$, which is also an $sr$-dimensional vector subspace of $\Fqm$ over $\Fp$. Set $R^*=R \sm \{ 0 \}$. By Equation~(\ref{eqn-T}), we have 
\begin{equation*}
R^*=\{ \ga^t \om \mid t \in T, \om \in \Fq^* \},
\end{equation*}
with $|R^*|=q^r-1$. By Equation~(\ref{eqn-DT3}), we have
\begin{equation}
  (\psi_a \times \vp_b)(D_{r,T}) = - \qml \sum_{u \in R^*} \zp^{\Trqmp(\ga^{f_b-f_a} u)} - \frac{(q^r-1)(\qml-1)}{q^m-1}-1       
\label{eqn-DT4}                                                                         
\end{equation}
Let 
$$
R^{\perp} = \{ x \in \Fqm \mid \Trqmp(xy) = 0 \mbox{ for each } y \in R \}
$$
be the dual space of $R$ with respect to the trace function $\Trqmp$, which is an $s(m-r)$-dimensional vector subspace of $\Fqm$ over $\Fp$. We note that
$$
\sum_{u \in R^*} \zp^{\Trqmp(\ga^{f_b-f_a} u)}=\begin{cases}
                                                                              q^r-1 & \mbox{if $\ga^{f_b-f_a} \in R^{\perp}$} \\
                                                                              -1      & \mbox{if $\ga^{f_b-f_a} \not\in R^{\perp}$}
                                                                          \end{cases}
$$
Moreover, recall that 
\begin{align*}
\ga^{f_a}=&\Nqmlqm(\al)^{f_a}=\Nqmlqm(\al^{f_a})=\Nqmlqm(a), \\ 
\ga^{f_b}=&\Nqmloqm(\be)^{f_b}=\Nqmloqm(\be^{f_b})=\Nqmloqm(b).
\end{align*}
Then we have
\begin{equation}
\sum_{u \in R^*}\zp^{\Trqmp(\ga^{f_b-f_a} u)} =\begin{cases}
                                                                              q^r-1 & \mbox{if $\ga^{f_b-f_a}=\frac{\Nqmloqm(b)}{\Nqmlqm(a)} \in R^{\perp}$} \\
                                                                              -1      & \mbox{if $\ga^{f_b-f_a}=\frac{\Nqmloqm(b)}{\Nqmlqm(a)} \not\in R^{\perp}$}
                                                                        \end{cases} 
\label{eqn-Rsum}                                                                        
\end{equation}
Combining Equations~(\ref{eqn-DT4}),~(\ref{eqn-Rsum}), we have 
\begin{align*}
  (\psi_a \times \vp_b)(D_{r,T}) =& - \qml \sum_{u \in R^*}\zp^{\Trqmp(\ga^{f_b-f_a} u)} - \frac{(q^r-1)(\qml-1)}{q^m-1}-1 \\
                                             =& \begin{cases}
                                                     -\frac{(q^r-1)(q^{m(\ell+1)}-1)}{q^m-1}-1 & \mbox{if $\frac{\Nqmloqm(b)}{\Nqmlqm(a)} \in R^{\perp}$,} \\
                                                     \frac{(q^m-q^r)(q^{m\ell}-1)}{q^m-1} & \mbox{if $\frac{\Nqmloqm(b)}{\Nqmlqm(a)} \not\in R^{\perp}$.}
                                                   \end{cases}                                 
\end{align*}

Combining Cases 1, 2, and 3, we conclude
\begin{align} 
    & (\psi_a \times \vp_b)(D_{r,T}) \nonumber \\
 = & \begin{cases} 
          -\frac{(q^r-1)(q^{m(\ell+1)}-1)}{q^m-1}-1 & \mbox{if $a \ne 0$, $b=0$ or $ab \ne 0$, $\frac{\Nqmloqm(b)}{\Nqmlqm(a)} \in R^{\perp}$,} \\
           \frac{(q^m-q^r)(q^{m\ell}-1)}{q^m-1} & \mbox{if $a=0$, $b \ne 0$ or $ab \ne 0$, $\frac{\Nqmloqm(b)}{\Nqmlqm(a)} \not\in R^{\perp}$.}
        \end{cases}
\label{eqn-DT5}         
\end{align}
Therefore, $D_{r,T}$ is a $(v,k_r,\la_r,\mu_r)$-PDS in $G$. 

Note that for each $(a,b) \in \bigcup_{i=0}^{e-1} (\Cqml_i \times \cup_{t \in T} \Cqmlo_{i+t})$,
$$
\frac{\Nqmloqm(b)}{\Nqmlqm(a)} \in \{\ga^t \om \mid t \in T, \om \in \Fq^* \}=R^*.
$$
The alternative expression (\ref{eqn-DTalt}) of $D_{r,T}$ then follows.

By Lemma~\ref{lem-connection}, $\Cay(G,D_{r,T})$ is a $(v,k_r,\la_r,\mu_r)$-SRG. Let $S$ be a subset of $G$ forming a clique in $\Cay(G,D_{r,T})$. In view of Remark~\ref{remark-Cayley} and Equation~(\ref{eqn-DTalt}), elements in $S$ must have distinct first coordinates. Therefore, $|S| \le q^{m\ell}$. On the other hand, $\Fqml \times \{ 0 \}$ corresponds to a clique of size $q^{m\ell}$ in $\Cay(G,D_{r,T})$. Thus, $\Cay(G,D_{r,T})$ contains a maximum clique $\Fqml \times \{ 0 \}$ of size $q^{m\ell}$. 
\end{proof}

Now we proceed to prove Theorem~\ref{theorem-dualmain2}, which implies Theorem~\ref{theorem-dualmain}.

\begin{proof}[Proof of Theorem~\ref{theorem-dualmain2}]
In view of Lemma~\ref{lem-dual} and Equation~(\ref{eqn-DT4}), the dual $D_{r,T}^+$ of $D_{r,T}$ can be expressed as 
$$
D_{r,T}^+=\{ \psi_a \times \vp_b \mid ab \ne 0, \frac{\Nqmloqm(b)}{\Nqmlqm(a)} \not\in R^{\perp} \} \cup \{ \psi_0 \times \vp_b \mid b \in \Fqmlo^* \} \subset \wh{G}.
$$
Recall that there exists a canonical group isomorphism $f: \wh{G} \rightarrow G$ such that $f(\psi_a \times \vp_b)=(a,b)$. Define $\wti{D_{r,T}^+}$ to be the image of $D_{r,T}^+$ under $f$, namely, $\wti{D_{r,T}^+}=f(D_{r,T}^+)$. Then 
\begin{equation*}
\wti{D_{r,T}^+}=\{ (a,b) \mid ab \ne 0, \frac{\Nqmloqm(b)}{\Nqmlqm(a)} \not\in R^{\perp} \} \cup \{ (0,b) \mid b \in \Fqmlo^* \} \subset G.
\end{equation*} 
Clearly, $D_{r,T}^+$ is a $(v^+,k_r^+,\la_r^+,\mu_r^+)$-PDS in $\wh{G}$ if and only if $\wti{D_{r,T}^+}$ is a $(v^+,k_r^+,\la_r^+,\mu_r^+)$-PDS in $G$. Therefore,
for the sake of simplicity, we identify $\wti{D_{r,T}^+}$ and $D_{r,T}^+$, so that $D_{r,T}^+$ is regarded as a subset of $G$ and
\begin{equation*}
D_{r,T}^+=\{ (a,b) \mid ab \ne 0, \frac{\Nqmloqm(b)}{\Nqmlqm(a)} \not\in R^{\perp} \} \cup \{ (0,b) \mid b \in \Fqmlo^* \}.
\end{equation*} 
Note that 
$$
\Fqm \sm R^{\perp}=\{ \ga^t\om \mid t \in {T^\perp}^c, \om \in \Fq^* \}=\cup_{t \in {T^\perp}^c} \Cqm_t.
$$
The above equation leads to the alternative expression (\ref{eqn-DTpalt}). The expression (\ref{eqn-DTp}) then follows from (\ref{eqn-DTpalt}). The parameters $(v^+,k_r^+,\la_r^+,\mu_r^+)$ follows from Lemma~\ref{lem-dual}.

By Lemma~\ref{lem-connection}, $\Cay(G,D_{r,T}^+)$ is a $(v^+,k_r^+,\la_r^+,\mu_r^+)$-SRG. Let $S$ be a subset of $G$ forming a clique in $\Cay(G,D_{r,T}^+)$. In view of Remark~\ref{remark-Cayley} and Equation~(\ref{eqn-DTpalt}), elements in $S$ must have distinct second coordinates. Therefore, $|S| \le q^{m(\ell+1)}$. On the other hand, $\{0\} \times \Fqmlo$ corresponds to a clique of size $q^{m(\ell+1)}$ in $\Cay(G,D_{r,T}^+)$. Thus, $\Cay(G,D_{r,T}^+)$ contains a maximum clique $\{0\} \times \Fqmlo$ of size $q^{m(\ell+1)}$. 
\end{proof}

\section{Concluding Remarks}
\label{sec-conclusion}

In this paper, we construct two classes of strongly regular Cayley graphs and PDSs based on elementary abelian groups. Respectively, they have the generalized Denniston and the dual generalized Denniston parameters, and hence are distinct from the well known Latin square type and negative Latin square type families. Our construction unifies and subsumes a number of existing constructions presented in the context of SRGs, PDSs, projective sets, and projective two-weight codes. Below, we mention several open problems which might be of interest for further research.

{\bf Problem 1}: As observed in Remarks~\ref{remark-main}, \ref{remark-code}, and \ref{remark-projset}, the main result Theorem~\ref{theorem-main} and its corollaries exactly recover some existing constructions as their special cases. On the other hand, for some previous constructions, we only know they share the same parameters as our construction. Therefore, a natural question is to study the whether these constructions with the same parameters are equivalent/isomorphic to each other. Specifically, 
\begin{itemize}
\item[(1)] Denniston's original construction uses an irreducible quadratic form from $\F_{2^m} \times \F_{2^m}$ to $\F_{2^m}$ \cite[Construction]{Denniston}. By identifying the additive group $(\F_{2^m} \times \F_{2^m}, +)$ with the additive group of the finite field $\F_{2^{2m}}$, we can employ the irreducible quadratic form $f : \F_{2^{2m}} \to \F_{2^m}$ given by $f(x) = x^{2^m+1}$ in Denniston's construction. This special case coincides with the construction in Theorem~\ref{theorem-main} with $p = 2$, $s = 1$, $\ell = 1$, and $1 \le r \le m-1$.
\item[(2)] It remains unclear whether the projective two-weight codes arising from Corollary~\ref{corollary-code}(1) with $\ell = 1$ are equivalent to the projective two-weight codes constructed by Bierbrauer and Edel \cite{BE}. One challenge in resolving this equivalence problem is the highly intricate nature of the construction in \cite{BE}.
\end{itemize}

{\bf Problem 2}: There has been intensive research on SRGs and PDSs with Denniston parameters in elementary abelian groups \cite{BXZ,BE,DHJP,Denniston,DeWinter,OttDM}. We note that SRGs and PDSs with Denniston parameters have also been constructed in nonelementary abelian groups \cite{DX} and nonabelian groups \cite{Brady,SDPS}. A natural question is whether SRGs and PDSs with generalized Denniston and dual generalized Denniston parameters can be constructed in groups other than elementary abelian groups.

{\bf Problem 3}: As observed in Remark~\ref{remark-relation}, SRGs and PDSs with generalized Denniston and dual generalized Denniston parameters form two rich families that are largely distinct from the Latin square type and the negative Latin square type families. It would be interesting to construct SRGs and PDSs belonging to new parameter families.

{\bf Problem 4}: In \cite{OttDCC}, Ott described a comprehensive framework involving a sophisticated treatment of Jacobi sums that generate PDSs in groups of the form $(\F_{q_1} \times \F_{q_2} \times \cdots \times \F_{q_s},+)$, where $q_i$, $1 \le i \le s$, are prime powers. Indeed, as mentioned in Remark~\ref{remark-main2}(3), \cite[Theorem 20]{OttDCC} provides a recursive construction of PDSs in group $(\F_{q_1} \times \F_{q_2},+)$, where $q_1$ and $q_2$ are not necessarily identical. Following this spirit, it is natural to consider adapting the approach in the current paper in order to construct PDSs in $(\F_{q_1} \times \F_{q_2},+)$ with $q_1 \ne q_2$. We note that there has been some interesting research along this direction \cite{FKM,OttDCC}.



\section*{Acknowledgements}

The authors would like to express their gratitude to Qing Xiang for enlightening comments on the blow up of point sets in $\PG(2,q^m)$, which lead to the paragraphs succeeding Remark~\ref{remark-projset}. The authors would like to sincerely thank the anonymous reviewers for their thoughtful and constructive comments. Shuxing Li's research was supported by the U.S. National Science Foundation Grant DMS-2452236 and the University of Delaware Research Foundation Strategic Initiative (UDRF-SI) program. The work of the Laura Johnson was supported by the Additional Funding Programme for Mathematical Sciences, delivered by EPSRC (EP/V521917/1) and the Heilbronn Institute for Mathematical Research.


%
%

\bibliographystyle{alphaurl}
\bibliography{reference}

@article{BBM,
  author  = {S. Ball and A. Blokhuis and F. Mazzocca},
  title   = {Maximal arcs in Desarguesian planes of odd order do not exist},
  journal = {Combinatorica},
  volume  = {17},
  number  = {1},
  pages   = {31--41},
  year    = {1997},
}

@article{BXZ,
  author  = {J. Bao and Q. Xiang and M. Zhao},
  title   = {Partial Difference Sets with Denniston parameters in elementary abelian $p$-groups},
  journal = {Finite Fields Appl.},
  volume  = {101},
  pages   = {Paper No. 102539, 14 pp},
  year    = {2025},
}

@book{BEW,
  author    = {B.C. Berndt and R.J. Evans and K.S. Williams},
  title     = {Gauss and Jacobi sums},
  series    = {Canadian Mathematical Society Series of Monographs and Advanced Texts},
  publisher = {John Wiley \& Sons, Inc.},
  address   = {New York},
  year      = {1998},
}

@article{BE,
  author  = {J. Bierbrauer and Y. Edel},
  title   = {A family of 2-weight codes related to BCH-codes},
  journal = {J. Combin. Des.},
  volume  = {5},
  number  = {5},
  pages   = {391--396},
  year    = {1997},
}

@article{Brady,
  author  = {A.C. Brady},
  title   = {Negative {L}atin square type partial difference sets in nonabelian groups of order 64},
  journal = {Finite Fields Appl.},
  volume  = {81},
  pages   = {Paper No. 102044, 11 pp},
  year    = {2022},
}

@book{BM,
  author    = {A.E. Brouwer and H. Van Maldeghem},
  title     = {Strongly regular graphs},
  series    = {Encyclopedia of Mathematics and its Applications},
  volume    = {182},
  publisher = {Cambridge University Press},
  address   = {Cambridge},
  year      = {2022},
}

@article{BWX,
  author  = {A.E. Brouwer and R.M. Wilson and Q. Xiang},
  title   = {Cyclotomy and strongly regular graphs},
  journal = {J. Algebraic Combin.},
  volume  = {10},
  number  = {1},
  pages   = {25--28},
  year    = {1999},
}

@article{CK,
  author  = {R. Calderbank and W.M. Kantor},
  title   = {The geometry of two-weight codes},
  journal = {Bull. London Math. Soc.},
  volume  = {18},
  number  = {2},
  pages   = {97--122},
  year    = {1986},
}

@article{Cameron,
  author  = {P.J. Cameron},
  title   = {Random strongly regular graphs?},
  journal = {Discrete Math.},
  volume  = {273},
  number  = {1-3},
  pages   = {103--114},
  year    = {2003},
}

@article{DHJP,
  author  = {J.A. Davis and S. Huczynska and L. Johnson and J. Polhill},
  title   = {Denniston partial difference sets exist in the odd prime case},
  journal = {Finite Fields Appl.},
  volume  = {99},
  pages   = {Paper No. 102499, 13 pp},
  year    = {2024},
}

@article{DX,
  author  = {J.A. Davis and Q. Xiang},
  title   = {A family of partial difference sets with Denniston parameters in nonelementary abelian 2-groups},
  journal = {European J. Combin.},
  volume  = {21},
  number  = {8},
  pages   = {981--988},
  year    = {2000},
}

@article{Denniston,
  author  = {R.H.F. Denniston},
  title   = {Some maximal arcs in finite projective planes},
  journal = {J. Combinatorial Theory},
  volume  = {6},
  pages   = {317--319},
  year    = {1969},
}

@article{DDHM,
  author  = {De Clerck, F. and Delanote, M. and Hamilton, N. and Mathon, R.},
  title   = {Perp-systems and partial geometries},
  journal = {Adv. Geom.},
  volume  = {2},
  number  = {1},
  pages   = {1--12},
  year    = {2002},
}

@article{DeWinter,
  author  = {De Winter, S.},
  title   = {Projective two-weight sets of Denniston type},
  journal = {Comb. Theory},
  volume  = {5},
  number  = {1},
  pages   = {Paper No. 2, 9 pp},
  year    = {2025},
}

@article{FKM,
  author  = {G.A. Fern\'{a}ndez-Alcober and R. Kwahira and L. Mart\'{\i}nez},
  title   = {Cyclotomy over products of finite fields and combinatorial applications},
  journal = {European J. Combin.},
  volume  = {31},
  number  = {6},
  pages   = {1520--1538},
  year    = {2010},
}

@book{Gauss,
  author    = {C.F. Gauss},
  title     = {Disquisitione Arithmeticae},
  publisher = {Yale University Press},
  address   = {New Haven, Conn.-London},
  year      = {1966},
  note      = {Translated into English by Arthur A. Clarke, S. J},
}

@article{JL21,
  author  = {J. Jedwab and S. Li},
  title   = {Packings of partial difference sets},
  journal = {Comb. Theory},
  volume  = {1},
  pages   = {Paper No. 18, 41 pp},
  year    = {2021},
}

@incollection{JL24,
  author    = {J. Jedwab and S. Li},
  title     = {Group rings and character sums: tricks of the trade},
  booktitle = {New advances in designs, codes and cryptography},
  pages     = {241--266},
  series    = {Fields Institute Communications},
  volume    = {86},
  publisher = {Springer},
  address   = {Cham},
  year      = {2024},
}

@book{LN,
  author    = {R. Lidl and H. Niederreiter},
  title     = {Finite fields},
  series    = {Encyclopedia of Mathematics and its Applications},
  volume    = {20},
  publisher = {Cambridge University Press},
  address   = {Cambridge},
  year      = {1997},
}

@article{Ma,
  author  = {S.L. Ma},
  title   = {A survey of partial difference sets},
  journal = {Des. Codes Cryptogr.},
  volume  = {4},
  number  = {3},
  pages   = {221--261},
  year    = {1994},
}

@article{M,
  author  = {K. Momihara},
  title   = {Certain strongly regular {C}ayley graphs on {$\Bbb{F}_{2^{2(2s+1)}}$} from cyclotomy},
  journal = {Finite Fields Appl.},
  volume  = {25},
  pages   = {280--292},
  year    = {2014},
}

@incollection{MWX,
  author    = {K. Momihara and Q. Wang and Q. Xiang},
  title     = {Cyclotomy, difference sets, sequences with low correlation, strongly regular graphs and related geometric substructures},
  booktitle = {Combinatorics and finite fields---difference sets, polynomials, pseudorandomness and applications},
  pages     = {173--198},
  series    = {Radon Series on Computational and Applied Mathematics},
  volume    = {23},
  publisher = {De Gruyter},
  address   = {Berlin},
  year      = {2019},
  note      = {Volumn 23 of Radon Series on Computational and Applied Mathematics},
}

@article{MX,
  author  = {K. Momihara and Q. Xiang},
  title   = {Lifting constructions of strongly regular Cayley graphs},
  journal = {Finite Fields Appl.},
  volume  = {26},
  pages   = {86--99},
  year    = {2014},
}

@article{OttDM,
  author  = {U. Ott},
  title   = {A generalization of a cyclotomic family of partial difference sets given by Fern\'{a}ndez-Alcober, Kwashira, and Mart\'{\i}nez},
  journal = {Discrete Math.},
  volume  = {339},
  number  = {8},
  pages   = {2153--2156},
  year    = {2016},
}

@article{OttDCC,
  author  = {U. Ott},
  title   = {On Jacobi sums, difference sets and partial difference sets in Galois domains},
  journal = {Des. Codes Cryptogr.},
  volume  = {80},
  number  = {2},
  pages   = {241--281},
  year    = {2016},
}

@article{P,
  author  = {J. Polhill},
  title   = {Generalizations of partial difference sets from cyclotomy to nonelementary abelian {$p$}-groups},
  journal = {Electron. J. Combin.},
  volume  = {15},
  number  = {1},
  pages   = {Research Paper 125, 13 pp},
  year    = {2008},
}

@book{Storer,
  author    = {T. Storer},
  title     = {Cyclotomy and difference sets},
  series    = {Lectures in Advanced Mathematics},
  volume    = {2},
  publisher = {Markham Publishing Co.},
  address   = {Chicago, IL},
  year      = {1967},
}

@misc{SDPS,
  author       = {E. Swartz and J.A. Davis and J. Polhill and K.W. Smith},
  title        = {Combinatorial transfer: a new method for constructing infinite families of nonabelian difference set, partial difference sets, and relative difference sets},
  howpublished = {arXiv:2407.18385},
  year         = {2024},
}

\end{document}